\documentclass[a4paper,12pt,final]{article}

\usepackage{amsmath, amssymb,amsthm}
\usepackage{mathrsfs}
\usepackage{enumerate}

\numberwithin{equation}{section}

\newcommand{\B}{\ensuremath{\mathbb{B}}}
\newcommand{\C}{\ensuremath{\mathbb{C}}}
\newcommand{\R}{\ensuremath{\mathbb{R}}}
\newcommand{\dv}{\mathop{\mathrm{div}}\nolimits}
\newcommand{\curl}{\mathop{\mathrm{curl}}\nolimits}
\newcommand{\supp}{\mathop{\mathrm{supp}}\limits}
\newcommand{\loc}{\text{\normalfont loc}}
\newcommand{\vv}{\boldsymbol{v}}
\newcommand{\vu}{\boldsymbol{u}}
\newcommand{\vw}{\boldsymbol{w}}
\newcommand{\vz}{\boldsymbol{z}}

\newcommand{\vB}{\boldsymbol{B}}

\newcommand{\vL}{\boldsymbol{L}}
\newcommand{\vW}{\boldsymbol{W}}
\newcommand{\vC}{\boldsymbol{C}}

\newcommand{\va}{\boldsymbol{a}}
\newcommand{\vb}{\boldsymbol{b}}
\newcommand{\vnu}{\boldsymbol{\nu}}
\newcommand{\vf}{\boldsymbol{f}}
\newcommand{\vg}{\boldsymbol{g}}
\newcommand{\vh}{\boldsymbol{h}}
\newcommand{\+}{|\!|\!|}%
\newcommand{\LN}{[\![}
\newcommand{\RN}{]\!]}

\theoremstyle{plain}
\newtheorem{theorem}{Theorem}[section]

\newtheorem{lemma}[theorem]{Lemma}
\newtheorem{proposition}[theorem]{Proposition}
\theoremstyle{definition}

\theoremstyle{remark}
\newtheorem{remark}[theorem]{Remark}

\begin{document}
\title{\bfseries\Large
On an existence theorem of global strong solution to the
magnetohydrodynamic system in three dimensional exterior domain}
\author{Norikazu Yamaguchi\\[.5cm]
{\normalsize Department of Mathematical Sciences,}\\
{\normalsize School of Science and Engineering, Waseda University}\\
{\normalsize 3-4-1 \=Okubo, Shinjuku-ku, Tokyo 169-8555, Japan}}
\date{}

\maketitle

\begin{abstract}
 In this paper we study the initial-boundary value problem for 
 the magnetohydrodynamic system in three dimensional exterior domain.
 We show an existence theorem of global in time strong solution
 for small $L^3$-initial data and we also show its asymptotic behavior when
 time goes to infinity.\\

 
\end{abstract}

\section{Introduction and main results}\label{sec:intro}
Let $\mathcal{O}$ be a \textit{simply connected} and bounded open set
in $\R^3$ with $C^{2,1}$-boundary. We choose some $R_0 > 0$
such that $\mathcal{O} \subset B_{R_0} = \{ x \in \R^3\,|\, |x| < R_0\}$
and fix it. 
Let $\varOmega$ be the exterior domain to $\mathcal{O}$, i.e., 
$\varOmega = \R^3 \setminus \overline{\mathcal{O}}$.
In this paper we are concerned with the initial-boundary value problem
of the magnetohydrodynamic system (the Ohm-Navier-Stokes system)
concerning the velocity $\vv=(v_1(x,t),v_2(x,t),v_3(x,t))$,
pressure $p=p(x,t)$ and 
magnetic field $\vB=(B_1(x,t),B_2(x,t),B_3(x,t))$ 
in $\varOmega \times (0,\infty)$:
\begin{equation}
\left\{
 \begin{aligned}
 &\vv_t  - \Delta \vv + (\vv \cdot \nabla) \vv + \nabla p + \vB \times \curl{\vB} =0 &\ &\text{in} &\ &\varOmega \times (0,\infty),\\
 &\vB_t + \curl{\curl{\vB}} + (\vv \cdot \nabla) \vB - (\vB \cdot \nabla)\vv =0 &\ &\text{in} &\ &\varOmega \times (0,\infty),\\
  &\dv{\vv} =0, \quad \dv{\vB} =0  &\ &\text{in} &\ &\varOmega \times
  (0,\infty),\\
  &\vv = 0, \quad  \vnu \cdot \vB =0, \quad \curl{\vB} \times \vnu =0,  &\ &\text{on} &\ &\partial\varOmega \times (0,\infty),\\
  &\vv(x,0) = \boldsymbol{a}, \quad \vB(x,0) = \boldsymbol{b}  &\ &\text{in} &\ &\varOmega.
 \end{aligned}
\right.
\tag{MHD}\label{eq:MHD}
\end{equation}
Here $\boldsymbol{a}=(a_1(x),a_2(x),a_3(x))$ and $\boldsymbol{b} =
(b_1(x)),b_2(x),b_3(x))$ are the prescribed initial data for the
velocity and magnetic field, respectively and $\vnu =
(\nu_1,\nu_2,\nu_3)$ is the unit outer normal on $\partial \varOmega$.
The magnetohydrodynamic system is known to be one of the mathematical models
describing the motion of the incompressible viscous and electrically
conducting Newtonian fluids.
This system is a coupled system of the Navier-Stokes system, Maxwell's
equations and Ohm's law under the MHD approximation (see e.g., Landau
and Lifshitz \cite{L-L}).

On the nonstationary problem of the magnetohydrodynamic system,  there
are many works when $\varOmega = \R^3$ or $\varOmega$ is bounded.
For example, Ladyzhenskaya and Solonnikov \cite{L-S-1960},
Duvaut and J.-L.~Lions \cite{D-L-1972} and Sermange and Temam \cite{S-T}.
However, all of the works above are done in the $L^2$ setting.
While on the other hand, Yoshida and Giga \cite{Y-G} studied
\eqref{eq:MHD} when $\Omega$ is bounded by analytic semigroup approach
similar to Giga and Miyakawa \cite{G-M-85} and they constructed the
unique global strong solution if the initial data $(\va,\vb)$ are
sufficiently small in sense of $L^3$.
In the exterior domain case, Kozono \cite{Kozono-87} showed the energy
decay of the weak solution of \eqref{eq:MHD}.
As far as the author knows, there has been no work on a global
in time existence of strong solution to \eqref{eq:MHD} when $\Omega$ is 
exterior domain.

For the nonstationary problem of the Navier-Stokes equations for the
motion of the viscous incompressible fluids, 
T.~Kato \cite{Kato-84} showed the global solvability of the Cauchy
problem if initial velocity $\va$ is sufficiently small with respect to
$L^n$-norm ($n \geq 2$ denotes the dimension).
The argument of Kato is based on the estimates of various $L^q$-norm of
the Stokes semigroup (in the whole space, the Stokes semigroup is
essentially the same as the heat semigroup $e^{t\Delta}$).
In particular, the $L^q$-$L^r$ type estimates for such semigroup play a
crucial role in his argument. 
The result of Kato was extended to the case of $n$-dimensional exterior
domain ($n \geq 3$) by Iwashita \cite{Iwashita}.
Iwashita showed the $L^q$-$L^r$ estimates for the Stokes semigroup in
exterior domain which will be introduced later and solved the initial
boundary value problem of the Navier-Stokes equations in exterior domain
by using Kato's iteration scheme. 
In view of Kato and Iwashita, if the initial value $(\va,\vb)$ are small
enough in the sense of the $L^3$-norm, we can expect that \eqref{eq:MHD}
admits a unique global strong solution.
Indeed, as mentioned before Yoshida and Giga \cite{Y-G} succeeded in
constructing the global $L^3$-solution when $\Omega$ is bounded domain.
Thus, our main purpose of the present paper is to show an existence
theorem of global strong solution for \eqref{eq:MHD}. 

Since the main point of the argument of Kato and Iwashita consists of
the study of the linearized problem.
Therefore in order to treat \eqref{eq:MHD} by such argument,
we have to study the linearized problems of \eqref{eq:MHD} and
 investigate the properties of solutions to such problem.
If we linearize \eqref{eq:MHD}, we obtain two systems of equations.
The first one is system of the Stokes equations and the second one is 
the following linear diffusion equations with the perfectly conducting
wall:
\begin{equation}
\left\{
	\begin{aligned}
		&\vu_t + \curl{\curl{\vu}} =0, \quad \dv{\vu} =0 &\, &\text{in}
	 &\, &\varOmega \times (0,\infty),\\
		&\vnu \cdot \vu = 0, \quad \curl{\vu} \times \vnu = 0 &\,
	 &\text{on} &\, &\partial\varOmega \times (0,\infty),\\
		&\vu(x,0)=\boldsymbol{b} &\, &\text{in} &\, &\varOmega.
	\end{aligned}
	\right.\label{eq:PW}
\end{equation}
For the nonstationary Stokes equations, we already had the $L^q$-$L^r$
estimates due to Iwashita, thus what we have to do here is to get
the $L^q$-$L^r$ estimates for the solutions of \eqref{eq:PW}.

To state main results of this paper precisely,
at this point we shall introduce notation used throughout this paper.
We use the following symbols for denoting the special sets,
$B_R = \{x \in \R^3\,|\, |x| < R\}$,
$S_R = \{x \in \R^3\,|\, |x| = R\}$,
$D_{L,R} = \{x \in \R^3\,|\, L \leq |x| \leq R\}$,
$\varOmega_R = \varOmega \cap B_R$,
$\partial \varOmega_R = \partial \varOmega \cup S_R$.

Let $D$ be any domain in $\R^3$.
For $1 \leq q \leq \infty$,
$L^q(D)$ denotes the usual Lebesgue space on $D$,
$W^{m,q}(D)$ denotes the usual $L^q$-Sobolev space of order $m$,
and $C_0^{\infty}(D)$ is the set of all infinitely differentiable functions
in $D$ with compact support in $D$.
For function spaces of vector valued functions,
we use the following symbols:
\begin{equation*}
 \vL^q(D)=\{ \vf = (f_1,f_2,f_3)\,|\, f_j \in L^q(D), j=1,2,3\},
\end{equation*}
likewise for $\vW^{m,q}(D)$, $\vC_0^{\infty}(D)$.
Moreover we define a function space $\vL^q_{R}(D)$ as follow:
\begin{align*}
 \vL^q_{R}(D) &= \{ \vf \in \vL^q(D)\,|\, \supp{\vf} \subset B_R\}.
\end{align*}
For the differentiation of three-vector of functions $\vf=(f_1,f_2,f_3)$ and
the scalar function $p$ we use the following symbols:
$\partial_j p = \partial p/\partial x_j$,
$p_t=\partial_t p=\partial p/\partial t$,
$\nabla p=(\partial_1 p,\partial_2 p,\partial_3 p)$,
\begin{gather*}
 \dv{\vf} = \sum_{j=1}^3 \partial_j f_j, \quad
 \curl{\vf} = (\partial_2 f_3 - \partial_3 f_2, \partial_3 f_1 - \partial_1 f_3,\partial_1 f_2 - \partial_2 f_1), \\
 \nabla^m \vf = (\partial^{\alpha}_{x} \vf \,|\, |\alpha|=m).
\end{gather*}
To denote various constants,
we use the same letters $C$ and $C_{A,B,\dots}$ means that the
constant depends on $A,B,\dots$.
The constants $C$ and $C_{A,B,\dots}$ may change from line to line.

In order to give an operator theoretic interpretation of \eqref{eq:MHD},
here we shall introduce the well known Helmholtz decomposition of
$\vL^q(\Omega)$. 
First, we shall introduce the following function space:
\begin{equation*}
C_{0,\sigma}^{\infty}(\varOmega) = \{ \vf \in
 \vC_{0}^{\infty}(\varOmega)\,|\, \dv{\vf} = 0 \text { in } \varOmega\}.
\end{equation*}
Let $1 < q < \infty$.
As is well known that the Banach space $\vL^q(\varOmega)$ admits
the Helmholtz decomposition (see Miyakawa \cite{My2}, Galdi
\cite[Chapter III]{GlI} and Simader and Sohr \cite{MR1190728}): 
\begin{equation*}
	\vL^q(\varOmega) = L^q_{\sigma}(\varOmega) \oplus G^q(\varOmega), \quad \oplus : \text{direct sum}.
\end{equation*}
Here 
\begin{align*}
	L^q_{\sigma}(\varOmega) &= \overline{C_{0,\sigma}^{\infty}(\varOmega)}^{\|\cdot\|_{{L^q}(\varOmega)}},\\
	G^q(\varOmega) &= \{\vf \in \vL^q(\varOmega)\,|\, \vf = \nabla p \text{
 for some } p \in L^q_{\loc}(\overline{\varOmega})\}.
\end{align*}
Since $\partial \varOmega$ is $C^{2,1}$-hypersurface,
the solenoidal space $L^q_{\sigma}(\varOmega)$ is characterized as (see
e.g., Galdi \cite{GlI}) 
\begin{equation}
 	L^q_{\sigma}(\varOmega) = \{\vf \in \vL^q(\varOmega)\,|\, \dv{\vf}=0
	 \text{ in } \varOmega,\ \vnu \cdot \vf =0 \text{ on } \partial\varOmega \}.
	 \label{eq:solenoidal-sp}
\end{equation}
Let $P = P_{q,\varOmega}$ be a continuous projection from $\vL^q(\varOmega)$
onto $L^q_{\sigma}(\varOmega)$ and then
\begin{equation}
 \|P \vf\|_{L^q(\varOmega)} \leq C_q \|\vf\|_{L^q(\varOmega)}
 \label{eq:HP-Boundedness}
\end{equation}
for any $\vf \in \vL^q(\varOmega)$.
Let us define the linear operators
$A = A_{q,\varOmega}$ and $\mathcal{M} = \mathcal{M}_{q,\varOmega}$ as
follows: 
\begin{equation*}
\begin{aligned}
 \mathcal{D}(A) &= L^q_{\sigma}(\varOmega) \cap \vW^{2,q}(\varOmega) \cap \vW^{1,q}_0(\varOmega),\\
 A \vv &= - P \Delta \vv \text{ for } \vv \in \mathcal{D}(A),\\
 \mathcal{D}(\mathcal{M}) &= L^q_{\sigma}(\varOmega) \cap \{\vB \in \vW^{2,q}(\varOmega) \, | \, \curl{\vB} \times \vnu = 0 \text{ on } \partial \varOmega\},\\
 \mathcal{M} \vB &= \curl{\curl{\vB}} \text{ for } \vB \in \mathcal{D}(\mathcal{M}).
\end{aligned}
\end{equation*}
The operator $A$ is usually called \textit{the Stokes operator} with non
slip boundary condition.
We note that the operator $\mathcal{M}$ is mapping 
from $\mathcal{D}(\mathcal{M})$ to $L^q_{\sigma}(\varOmega)$.
By using $A$ and $\mathcal{M}$, \eqref{eq:MHD} is rewritten by the
following Cauchy problem of abstract evolution equations in the Banach
space $L^q_{\sigma}(\Omega) \times L^q_{\sigma}(\Omega)$:
\begin{equation}
 \left\{
 \begin{aligned}
  &\frac{d\vv(t)}{dt} + A \vv(t) +P[(\vv(t) \cdot \nabla)\vv(t) - (\vB(t)
  \cdot \nabla) \vB(t)] =0, &\quad &t>0,\\
  &\frac{d\vB(t)}{dt} + \mathcal{M} \vB(t) + (\vv(t) \cdot \nabla)\vB(t)
  - (\vB(t) \cdot \nabla) \vv(t) = 0, & &t>0,\\
  &\vv(0) = \va, \quad \vB(0) = \vb.
 \end{aligned}
 \right.
 \tag{ACP}\label{eq:ACP}
\end{equation}
Here we have used the well known  formula:
\begin{equation*}
 \vB \times \curl{\vB} = - (\vB \cdot \nabla) \vB + \frac{\nabla |\vB|^2}{2}.
\end{equation*}
The second term in the right hand side of the above relation
is eliminated by the Helmholtz projection $P$.
According to Miyakawa \cite{My2} and Borchers and Sohr \cite{B-S},
$-A$ generates a bounded analytic semigroup $(e^{-tA})_{t \geq 0}$ on
$L^q_{\sigma}(\varOmega)$ and according to Miyakawa \cite{My} and 
Shibata and Yamaguchi \cite{MHD-local-energy} the 
operator $-\mathcal{M}$ also generates a bounded analytic semigroup
$(e^{-t \mathcal{M}})_{t \geq 0}$ on $L^q_{\sigma}(\varOmega)$.
Therefore, by virtue of Duhamel's principle, \eqref{eq:ACP} is converted
into the following system of integral equations:
\begin{equation}
\left\{
 \begin{aligned}
  \vv(t) &= e^{-t A}\va - \int_{0}^{t} e^{-(t-s)A} P[(\vv(s) \cdot \nabla)\vv(s) - (\vB(s) \cdot \nabla) \vB(s)]\,ds,\\
  \vB(t) &= e^{-t \mathcal{M}}\vb - \int_{0}^{t} e^{-(t-s)\mathcal{M}}
  [(\vv(s) \cdot \nabla)\vB(s) - (\vB(s) \cdot \nabla) \vv(s)]\,ds. 
 \end{aligned}
\right.\tag{INT}
\label{eq:INT}
\end{equation}
For notational simplicity,
we set 
$\vv_0(t) = e^{-t A}\va$, $\vB_0(t) = e^{-t \mathcal{M}}\vb$,
\begin{align*}
 F[\vv,\vB](t) &= - \int_{0}^{t} e^{-(t-s)A} P[(\vv(s) \cdot \nabla)\vv(s) - (\vB(s) \cdot \nabla) \vB(s)]\,ds,\\
 G[\vv,\vB](t) &= - \int_{0}^{t} e^{-(t-s)\mathcal{M}} [(\vv(s) \cdot \nabla)\vB(s) - (\vB(s) \cdot \nabla) \vv(s)]\,ds.
\end{align*}
Our aim of this paper is deduced to solve \eqref{eq:INT}
by contraction mapping principle (or Kato's iteration scheme).
In order to do this, we need $L^q$-$L^r$ estimates for 
the semigroups $e^{-t A}$ and $e^{-t \mathcal{M}}$.

\bigskip

We are now in a position to state our main results.
The first result is concerning $L^q$-$L^r$ estimates for the 
semigroup $e^{-t \mathcal{M}}$.
\begin{theorem}[$L^q$-$L^r$ estimates]
	\label{thm:lq-lr}
  \ \\
 {\normalfont (i)}
 Let $1 \leq q \leq r \leq \infty$ and $q \not=\infty, r\not=1$.
 Then there exists a constant $C=C_{q,r}>0$ such that
 \begin{equation*}
  \|e^{-t \mathcal{M}} \vf\|_{L^r(\varOmega)} \leq C
  t^{-\frac{3}{2}\left(\frac{1}{q}-\frac{1}{r}\right)}
  \| \vf \|_{L^q(\varOmega)},
  \quad t>0
 \end{equation*}
 for any $\vf \in L^q_{\sigma}(\varOmega)$.\\
 {\normalfont (ii)} 
 Let $1 \leq  q \leq r \leq 3$, $r \not=1$.
 Then there exists a constant $C=C_{q,r}>0$ such that
 \begin{equation*}
  \|\nabla e^{-t \mathcal{M}} \vf\|_{L^r(\varOmega)} \leq C
  t^{-\frac{3}{2}\left(\frac{1}{q}-\frac{1}{r}\right) - \frac{1}{2}}
  \|\vf\|_{L^q(\varOmega)},
  \quad t>0
 \end{equation*}
 for any $\vf \in L^q_{\sigma}(\varOmega)$.
\end{theorem}

The basic idea to prove Theorem~\ref{thm:lq-lr} is similar to that of
Iwashita \cite{Iwashita} for the Stokes semigroup.
Iwashita's idea is based on the local energy decay property of the
semigroup near the obstacle $\mathcal{O}$.
Such local energy decay estimate for $e^{-t\mathcal{M}}$ is obtained by 
Shibata and Yamaguchi \cite{MHD-local-energy} (see also \cite{5th-RIMS}).
\begin{theorem}
 [local energy decay \cite{MHD-local-energy}]
 \label{thm:local-energy}
 Let $1 < q < \infty$. For any $R > R_0$ ,
 there exists a constant $C=C_{q,R} > 0$ such that
 \begin{equation*}
  \|e^{-t \mathcal{M}} \vf\|_{W^{2,q}(\varOmega_R)} \leq C t^{-\frac{3}{2}} \| \vf \|_{L^q(\varOmega)}, \quad t \geq 1,
 \end{equation*}
 for any $ \vf \in L^q_{\sigma}(\varOmega) \cap \vL^q_{R}(\varOmega)$.
\end{theorem}

The following theorem by Iwashita \cite{Iwashita} is concerning the
$L^q$-$L^r$ estimates for the Stokes semigroup, which is refined by
Maremonti and Solonnikov \cite{M-S} and Enomoto and Shibata
\cite{Enomoto-Shibata-05-JMFM} (see also Giga and Sohr
\cite{MR991022}).
\begin{theorem}[$L^q$-$L^r$ estimates for the Stoke semigroup
 \cite{Enomoto-Shibata-05-JMFM,MR991022,Iwashita,M-S}] 
	\label{thm:lq-lr-Stokes}
  \ \\
 {\normalfont (i)}
 Let $1 \leq q \leq r \leq \infty$ and $q \not=\infty, r\not=1$.
 Then there exists a constant $C=C_{q,r}>0$ such that
 \begin{equation*}
  \|e^{-t A} \vf\|_{L^r(\varOmega)} \leq C
  t^{-\frac{3}{2}\left(\frac{1}{q}-\frac{1}{r}\right)}
  \|\vf\|_{L^q(\varOmega)},
  \quad t>0
 \end{equation*}
 for any $\vf \in L^q_{\sigma}(\varOmega)$.\\
 {\normalfont (ii)} 
 Let $1 < q \leq r \leq 3$.
 Then there exists a constant $C=C(q,r)>0$ such that
 \begin{equation*}
  \|\nabla e^{-t A} \vf\|_{L^r(\varOmega)} \leq C
  t^{-\frac{3}{2}\left(\frac{1}{q}-\frac{1}{r}\right) - \frac{1}{2}}
  \|\vf\|_{L^q(\varOmega)},
  \quad t>0
 \end{equation*}
 for any $\vf \in L^q_{\sigma}(\varOmega)$.
\end{theorem}

Finally, applying  Theorem~\ref{thm:lq-lr} and Theorem~\ref{thm:lq-lr-Stokes}
we obtain an existence theorem of global in time strong solution for
\eqref{eq:MHD} with small initial data.
\begin{theorem}[Global existence]
	\label{thm:global}
 There exists an $\eta=\eta(\varOmega)>0$ such that 
 if $(\va,\vb)\in L^3_{\sigma}(\varOmega) \times L^3_{\sigma}(\varOmega)$ satisfies
 $\|(\va,\vb)\|_{3} \leq \eta$ then \eqref{eq:MHD} has a unique global
 strong solution $(\vv(t),\vB(t)) \in BC([0,\infty);L^3_{\sigma}(\varOmega)
 \times L^3_{\sigma}(\varOmega))$ which possesses the followings{\normalfont :}
 \begin{align}
  &\lim_{t \rightarrow 0}\|(\vv(t),\vB(t) -
  (\va,\vb))\|_{L^3(\varOmega)} =0,\notag \\
  &\begin{aligned}
  &\lim_{t \rightarrow
  +0}t^{\frac{1}{2}-\frac{3}{2q}}\|(\vv(t),\vB(t))\|_{L^{q}(\varOmega)}\\
   &\ + \lim_{t \rightarrow +0} t^{\frac{1}{2}}\|\nabla
  (\vv(t),\vB(t))\|_{L^3(\varOmega)} =0 \quad\text{for\ } 3 < q < \infty;
   \end{aligned} \notag \\
  &\|(\vv(t),\vB(t))\|_{L^q(\Omega)} = o\left(t^{-\frac{1}{2}+\frac{3}{2q}}\right)\quad \text{for } 3 \leq q \leq \infty,\label{eq:AS-1}\\
  &\|\nabla(\vv(t),\vB(t))\|_{L^3(\Omega)} =
  o\left(t^{-\frac{1}{2}}\right)\label{eq:AS-2}. 
 \end{align}
 as $t \rightarrow \infty$.
 Here $BC(I;X)$ denotes the class of $X$-valued bounded and continuous
 function on interval $I$.
\end{theorem}
\begin{remark}
 We do not require any smallness assumption on the initial data
 for  proving  the local in time existence of solution to \eqref{eq:MHD}.
\end{remark}

Below, 
in section~\ref{sec:preliminaries} 
we prepare the well known Bogovski\u{\i}'s lemma and some lemmas which
will be used in the latter sections.
In section~\ref{sec:lq-lr} we shall prove Theorem~\ref{thm:lq-lr} with
aid of
$L^q$-$L^r$ estimates for the heat kernel,
Theorem~\ref{thm:local-energy} and cut-off technique. 
By using Theorems~\ref{thm:lq-lr} and \ref{thm:lq-lr-Stokes},
we prove Theorem~\ref{thm:global} in section~\ref{sec:global-existence}.

\section{Preliminaries}\label{sec:preliminaries}
In this section, we prepare some useful lemmas which will be used in 
the latter sections.
In Section~\ref{sec:lq-lr} we will prove Theorem~\ref{thm:lq-lr} by
cut-off technique.
In order to keep the divergence free condition in cut-off procedure,
we are due to the well known lemma by Bogovski\u{\i} \cite{Bg}
(see also Galdi \cite[Chapter~III]{GlI}).
In order to state Bogovski\u{\i}'s lemma, 
we shall introduce the function spaces $\dot{W}{}^{m,q}(D)$ and 
$\dot{W}{}^{m,q}_a(D)$ as follows:
\begin{gather*}
 \dot{W}{}^{m,q}(D) = \overline{C_0^{\infty}(D)}^{\| \cdot \|_{W^{m,q}}},\\
 \dot{W}{}^{m,q}_a (D) = \left\{f \in \dot{W}{}^{m,q}(D) \,\bigg|\,\int_{D} f(x)\,dx =
  0 \right\}.
\end{gather*}
Here $D$ stands for a bounded domain in $\R^3$ with smooth boundary
$\partial D$. 
We note that $\dot{W}{}^{0,q}(D)=L^q(D)$.
\begin{lemma}
 \label{lem:Bog}
 Let $1<q<\infty$ and let $m$ be a non-negative integer.
 Then there exists a bounded linear operator
 $\B \equiv \B_{D}: \dot{W}^{m,q}_a(D)  \rightarrow \dot{\vW}^{m+1,q} (\R^3)$ such that
 \begin{gather*}
  \supp{\B[f]} \subset D,\\
  \dv{\B[f]} = f \text{ in } \R^3.
 \end{gather*}
\end{lemma}
To use Lemma~\ref{lem:Bog}, we shall rely on the following lemma.
\begin{lemma}
 \label{prop:d-ext}
 Let $1<q<\infty$, $R > L > R_0$ and let $\varphi(x) \in C_0^{\infty}(\R^3)$
 such that $\varphi(x) = 1$ for $|x| \leq L$ and 
 $\varphi(x) = 0$ for $|x| \geq R$.
 \begin{enumerate}[{\normalfont (i)}]
  \item If $\vu \in \vW^{2,q}(\R^3)$ and $\vu$ satisfies the condition{\normalfont :} $\dv{\vu}=0$ in $\R^3$,
  then $(\nabla \varphi) \cdot \vu \in \dot{W}{}^{2,q}_{a}(D_{L,R})$.
  \item If $\vu \in \vW^{2,q}(\varOmega)$ and $\vu$ satisfies the conditions{\normalfont :}
   $\dv{\vu} = 0$ in $\varOmega$ and $\vnu \cdot \vu=0$ on $\partial\varOmega$,
   then $(\nabla \varphi) \cdot \vu \in \dot{W}{}^{2,q}_{a}(D_{L,R})$.
 \end{enumerate}
\end{lemma}

Next, we shall introduce the results in the case of bounded domain $D$.
From Akiyama, Kasai, Shibata and Tsutsumi \cite{AKST-04},
it follows the following proposition.
\begin{proposition}
 \label{prop:BD}
Let $1 < q < \infty$.
Assume that $\partial D \in C^{2,1}$.
Then for any $\vf \in \vL^q(D)$ %
there exists a unique solution $\vu \in \vW^{2,q}(D)$ of the following
 system{\normalfont:} 
\begin{equation*}
\left\{
 \begin{aligned}
  \vu - \Delta \vu &= \vf &\quad &\text{in} &\ &D,\\
  \curl{\vu} \times \vnu &=0 & &\text{on} & &\partial D,\\
   \vnu \cdot \vu &=0 & &\text{on} & &\partial D,
 \end{aligned}
\right.
\end{equation*}
 which satisfies the estimate{\normalfont :}
\begin{equation*}
  \| \vu \|_{W^{2,q}(D)} \leq C \| \vf \|_{L^q(D)}.
\end{equation*}
\end{proposition}
 
Next we shall introduce the resolvent estimate.
The resolvent problem corresponding to \eqref{eq:PW} is given by the
following Laplace system: 
\begin{equation}
\left\{
	\begin{aligned}
		\lambda \vu - \Delta \vu &=\vf &\, &\text{in} &\, &\varOmega,\\
		\curl{\vu} \times \vnu &=0 &\, &\text{on} &\, &\partial \varOmega,\\
		\vnu \cdot \vu &=0 &\, &\text{on} &\, &\partial \varOmega,
	\end{aligned}
\right.\label{eq:RP}
\end{equation}
The following theorem obtained by Akiyama, Kasai, Shibata and Tsutsumi
\cite{AKST-04} is concerned with the resolvent estimate for \eqref{eq:RP}.
\begin{theorem}
 \label{thm:resolvent-estimate}
 Let $1 < q < \infty$, $0 < \epsilon < \pi/2$ and $\delta > 0$.
 Set
 \begin{equation*}
  \Sigma_{\epsilon,\delta} = \{\lambda \in \C \setminus \{0\}\,|\, |\arg{\lambda}| \leq \pi - \epsilon, |\lambda| \geq \delta\}.
 \end{equation*}
  Then, for any $\vf \in L^q_{\sigma}(\varOmega)$ and $\lambda \in \Sigma_{\epsilon,\delta}$, \eqref{eq:RP}
   admits a unique solution $\vu \in \vW^{2,q}(\varOmega)$ possessing
 the estimate\,{\normalfont:}
 \begin{equation}
  |\lambda| \|\vu\|_{L^q(\varOmega)} + \|\vu\|_{W^{2,q}(\varOmega)} \leq C_{\epsilon,\delta}\|\vf\|_{L^q(\varOmega)}.
  \label{eq:Re-E}
 \end{equation}
\end{theorem}

On the linear operator $\mathcal{M}_q$ defined in
 Section~\ref{sec:intro},  we quote the following theorem due to 
Shibata and Yamaguchi \cite{MHD-local-energy}.
\begin{theorem}
 \label{thm:B-B*}
 Let $1 < q < \infty$, $q'=q/(q-1)$
 and $\mathcal{M}_q^{\ast}$ be an adjoint operator of $\mathcal{M}_q$.
 Then we have $\mathcal{M}_q^{\ast} = \mathcal{M}_{q'}$.
\end{theorem}

\section{Proof of Theorem~\ref{thm:lq-lr}}
\label{sec:lq-lr}
In this section we shall prove Theorem~\ref{thm:lq-lr}.
Our proof is based on the ideas due to Iwashita \cite{Iwashita} and
Hishida \cite{MR2085848}.
Here and hereafter $T(t)$ denotes the analytic semigroup generated
by $-\mathcal{M}_q$, i.e., $T(t) \equiv e^{-t \mathcal{M}}$.
Given $\vf \in L^q_{\sigma}(\varOmega)$, we set $\vu(t) = T(t)\vf$.
Then $\vu(t)$ solves the following initial-boundary value problem:
\begin{equation}
\left\{
	\begin{aligned}
		&\vu_t - \Delta \vu =0, \quad \dv{\vu} =0 &\, &\text{in}
	 &\, &\varOmega \times (0,\infty),\\
		&\vnu \cdot \vu = 0, \quad \curl{\vu} \times \vnu = 0 &\,
	 &\text{on} &\, &\partial\varOmega \times (0,\infty),\\
		&\vu(x,0) =\boldsymbol{\vf} &\, &\text{in} &\, &\varOmega.
	\end{aligned}
	\right.\label{eq:P-1}
\end{equation}
Here we have used the well known formula:
\begin{equation}
 \Delta \vu = \nabla \dv{\vu} - \curl{\curl{\vu}}.\label{eq:P-2}
\end{equation}

\subsection*{1st step}
As a first step, we shall show the following lemma.
\begin{lemma}
 \label{lem:P-1}
 Let $1 < q < \infty$ and $R > R_0 + 3$.
 Then there exists a $C=C_{q,\varOmega,R}>0$ such that
 \begin{equation*}
  \|\partial_t T(t) \vf\|_{W^{1,q}(\varOmega_R)} + \|T(t) \vf\|_{W^{2,q}(\varOmega_R)} \leq C t^{-\frac{3}{2q}} \|\vf\|_{L^q(\varOmega)}
 \end{equation*}
 for any $t \geq 2$ and $\vf \in L^q_{\sigma}(\varOmega)$.
\end{lemma}
\begin{proof}
 Since we consider the case when $t \geq 2$, we set
 \begin{equation}
  \vg = T(1)\vf, \quad \vv(t) = T(t)\vg = T(t+1)\vf.
  \label{eq:P-3}
 \end{equation}
 By \eqref{eq:Re-E}, the analytic semigroup theory (see e.g., Pazy
 \cite{Pazy}) and  \eqref{eq:P-1}, we have
 \begin{align}
  &\|\vg\|_{W^{2,q}(\varOmega)} \leq C \|\vf\|_{L^q(\varOmega)}, \quad
  \vg \in \mathcal{D}(\mathcal{M}),\label{eq:P-4}\\ 
  &\vv(t) \in C([0,\infty);\vW^{2,q}(\varOmega)) \cap
  C^1((0,\infty);\vL^q(\varOmega)),\label{eq:P-5}\\ 
  &\left\{\begin{aligned}
		&\vv_t - \Delta \vv =0, \quad \dv{\vv}=0 &\, &\text{in}
	 &\, &\varOmega \times (0,\infty),\\
		&\vnu \cdot \vv = 0, \quad \curl{\vv} \times \vnu = 0 &\,
	 &\text{on} &\, &\partial\varOmega \times (0,\infty),\\
		&\vv(x,0)=\vg &\, &\text{in} &\, &\varOmega.
	\end{aligned}
	\right.\label{eq:P-6}
 \end{align}
 Let $\psi \in C^{\infty}(\R^3)$ such that $\psi(x)=1$ for $|x| \geq
 R+1$  and $\psi(x) = 0$ for $|x| \leq R$.
 By \eqref{eq:P-4} and Lemma~\ref{prop:d-ext},
 we have  $(\nabla \psi) \cdot \vg \in W^{2,q}_{a}(D_{R,R+1})$ and 
 therefore by Lemma~\ref{lem:Bog} we have
 \begin{gather}
  \B_{D_{R,R+1}}[(\nabla \psi) \cdot \vg] \in W^{3,q}(\R^3), \quad
  \supp{\B_{D_{R,R+1}}[(\nabla \psi) \cdot \vg]} \subset D_{R,R+1},\notag \\
  \dv{\B_{D_{R,R+1}}[(\nabla \psi) \cdot \vg]} = (\nabla \psi) \cdot \vg,\notag \\
  \|\B_{D_{R,R+1}}[(\nabla \psi) \cdot \vg]\|_{W^{3,q}(\R^3)} 
    \leq C \|\vf\|_{L^q(\varOmega)}.\label{eq:P-7}
 \end{gather}
 In what follows, for notational simplicity,
 we use the abbreviation $\B = \B_{D_{R,R+1}}$.

 Let $E(t)$ be the Gaussian kernel, namely,  
 \begin{equation}
   E(t)=E(x,t) = \frac{1}{(4\pi t)^{\frac{3}{2}}} \exp\left(-\frac{|x|^2}{4t}\right)\label{eq:Gauss}
 \end{equation}
 and set 
 \begin{equation*}
  \vh = \psi \vg - \B[(\nabla \psi) \cdot \vg], \quad
  \vw = E(t) \ast \vh = \frac{1}{(4\pi t)^{3/2}} \int_{\R^3} \exp{\left(-\frac{|x-y|^2}{4t}\right)} \vh(y)  \,dy.
 \end{equation*}
 By \eqref{eq:P-4} and \eqref{eq:P-7}, we see that
 \begin{equation}
 \begin{gathered}
  \vh \in \vW^{2,q}(\R^3),\\
  \dv{\vh} = 0 \ \text{in}\ \R^3,\\
  \vh = \vg,\quad |x| \geq R+1 ,\\
  \|\vh\|_{W^{2,q}(\R^3)} \leq C_q \|\vf\|_{L^q(\varOmega)}.
 \end{gathered}
 \label{eq:P-8}
 \end{equation}
 Applying Young's inequality to $\vw(t)$ and using \eqref{eq:P-8}, we obtain
 \begin{align}
  &\vw(t) \in C([0,\infty);\vW^{2,q}(\R^3)) \cap C^1([0,\infty);\vL^q(\R^3)),
  \label{eq:P-9}\\
  &\vw_t - \Delta \vw =0, \quad \dv{\vw}=0 \quad \text{in}\ \R^3 \times (0,\infty),\quad \vw(0)=\vh,
  \label{eq:P-10}\\
  &\|\nabla^j \vw(t)\|_{L^r(\R^3)} \leq C_{q,r}
  (1+t)^{-\frac{3}{2}\left(\frac{1}{q}-\frac{1}{r}\right) - \frac{j}{2}}
  \|\vf\|_{L^q(\varOmega)},\quad j=1,2,\ t \geq  1,\notag \\
	&\|\vw_t\|_{L^r(\R^3)} + \|\nabla^2 \vw(t)\|_{L^r(\R^3)} \leq 
	 C_{q,r} (1+t)^{-\frac{3}{2}\left(\frac{1}{q}-\frac{1}{r}\right)-1} \|\vf\|_{L^q(\varOmega)}
	 \label{eq:P-11}
 \end{align}
 provided that $1 < q \leq r \leq \infty$.
 Since $\dv{\vw}=0$,
 by Lemma~\ref{lem:Bog} we have $(\nabla \psi) \cdot \vw(t) \in C([0,\infty);\dot{W}^{2,q}_a(D_{R,R+1}))$,
 and therefore we set
 \begin{equation}
  \vz(t) = \vv(t) - \psi \vw(t) + \B[(\nabla \psi) \cdot \vw(t)].
  \label{eq:P-12}
 \end{equation}
 Then, from \eqref{eq:P-5} and \eqref{eq:P-9} and Lemma~\ref{lem:Bog}
 we obtain
  \begin{align}
  &\vz(t) \in C([0,\infty);\vW^{2,q}(\varOmega)) \cap
   C^1([0,\infty);\vL^q(\varOmega)),\label{eq:P-13}\\ 
  &\left\{\begin{aligned}
		&\vz_t - \Delta \vz =\boldsymbol{F}(t), \quad \dv{\vz} =0 &\,
		   &\text{in} 
	 &\, &\varOmega \times (0,\infty),\\
		&\vnu \cdot \vz = 0, \quad \curl{\vz} \times \vnu = 0 &\,
	 &\text{on} &\, &\partial\varOmega \times (0,\infty),\\
		&\vz(0) =\vz_0 &\, &\text{in} &\, &\varOmega,
	\end{aligned}
	\right.\label{eq:P-14}
 \end{align}
 where we have set 
 \begin{equation}
 \begin{gathered}
  \boldsymbol{F}(t) = 2 \nabla \vw(t) \cdot \nabla \psi + (\Delta \psi) \vw(t) + (\partial_t - \Delta) \B[(\nabla \psi) \cdot \vw(t)],\\
   \vz_0 = \vg - \psi \vh + \B[(\nabla \psi) \cdot \vh].
 \end{gathered}
 \label{eq:P-15}
 \end{equation}
 We shall show that 
 \begin{gather}
 	\boldsymbol{F}(t) \in C([0,\infty);L^q_{\sigma}(\varOmega)), \quad 
 	 \supp{\boldsymbol{F}(t)} \subset D_{R,R+1} \quad \text{for any } t>0,
 	 \label{eq:P-16}\\
 	 \|\boldsymbol{F}(t)\|_{L^q(\varOmega)} \leq C(1+t)^{-\frac{3}{2q}} \|\vf\|_{L^q(\varOmega)},
 	 \label{eq:P-17}\\
 	 \vz_0 \in \mathcal{D}(\mathcal{M}_q), \quad \vz_0 = 0 \quad \text{for } x \not\in B_{R+1},
 	 \label{eq:P-18}\\
 	 \|\vz_{0}\|_{W^{2,q}(\varOmega)} \leq C \| \vf \|_{L^q(\varOmega)}.
 	\label{eq:P-19}
 \end{gather}
 In fact, since 
 \begin{equation*}
  (\partial_{t} - \Delta)(\psi \vw(t)) = - 2 \nabla \vw (t) \cdot \nabla \psi - (\Delta \psi) \vw(t),
 \end{equation*}
 by Lemma~\ref{prop:d-ext} we have
 \begin{align}
  \dv{\boldsymbol{F}(t)} &= - \dv{\{(\partial_{t} - \Delta)(\psi \vw(t))\}} + (\partial_{t} - \Delta) \dv{\B[(\nabla \psi) \cdot \vw(t)]} \notag \\
  &= - (\partial_{t} - \Delta)[\dv{(\psi \vw(t))} - (\nabla \psi) \cdot \vw(t)] = 0 \label{eq:P-20},
 \end{align}
 because $\dv{\vw(t)}=0$.
 Obviously, $\supp{\boldsymbol{F}(t)} \subset D_{R,R+1}$.
 In particular, we have $\vnu \cdot \boldsymbol{F}(t) =0$ on
 $\partial\varOmega$  for any $t \geq 0$,
 which combined with \eqref{eq:solenoidal-sp} and \eqref{eq:P-20} implies that
 $\boldsymbol{F}(t) \in L^q_{\sigma}(\varOmega)$ for any $t \geq 0$.
 Clearly, by \eqref{eq:P-5} and \eqref{eq:P-13},
 $\boldsymbol{F}(t) \in C([0,\infty);\vL^q(\varOmega))$,
 which completes the proof of \eqref{eq:P-16}.
 By Lemma~\ref{lem:Bog} and \eqref{eq:P-11} with $r = \infty$,
 we have
 \begin{align*}
  \|\boldsymbol{F}(t)\|_{L^q(\varOmega)}
  &\leq C_q 
    \{
		  \| |\nabla\psi| \nabla \vw(t)\|_{L^q(\varOmega)} + 
		  \||\Delta \psi| \vw(t)\|_{L^q(\varOmega)} \\
	&\qquad 
	   + \|\nabla \psi \cdot \vw(t)\|_{W^{1,q}(\varOmega)} + 
		  \|\nabla \psi \cdot \vw_t(t)\|_{L^q(\varOmega)}
		\}\\
  &\leq C_{q,R} \{ \|\vw(t)\|_{W^{1,\infty}(\R^3)} + \|\vw_t(t)\|_{L^{\infty}(\R^3)}  \}\\
  &\leq C_{q,R} (1+t)^{-\frac{3}{2q}} \|\vf\|_{L^q(\varOmega)}.
 \end{align*}
 By \eqref{eq:P-8} we see that $\vg = \psi \vh$ for $x \not \in B_{R+1}$.
 Furthermore, $\supp\B[(\nabla \psi) \cdot \vh] \subset D_{R,R+1}$.
 Therefore, by \eqref{eq:P-8} we have $\dv{\vz_0} = 0$ in $\varOmega$,
 $\|\vz_0\|_{W^{2,q}(\varOmega)} \leq C_q \|\vf\|_{L^q(\varOmega)}$ and
 $\vz_0 = 0$ for $x \not\in B_{R+1}$.
 Since $\vz_{0} = \vg$ for $|x| \leq R$,
 $\vg = T(1)\vf$ implies that $\vnu \cdot \vz_0 =0 $ and
 $\curl{\vz_0} \times \vnu = 0$ on $\partial\varOmega$.
 These facts imply that $\vz_0 \in \mathcal{D}(\mathcal{M}_q)$.
 Therefore we get \eqref{eq:P-16}, \eqref{eq:P-17} and \eqref{eq:P-18}.
 
 By \eqref{eq:P-13}, \eqref{eq:P-14} and \eqref{eq:P-16} and Duhamel's
 principle,  we have
 \begin{equation}
  \vz(t) = T(t)\vz_0 + \int_0^t T(t-s) \boldsymbol{F}(s)\,ds.
  \label{eq:P-21}
 \end{equation}
 Let $t \geq 1$.
 In view of \eqref{eq:P-16} and \eqref{eq:P-18}, we can apply
 Theorem~\ref{thm:local-energy} (local energy decay)  to estimate
 $\vz(t)$, and then we have 
 \begin{equation}
  \begin{aligned}
    \| \vz \|_{W^{1,q}(\varOmega_R)}
    &\leq C_{R} t^{-\frac{3}{2}} \|\vz_0\|_{L^q(\varOmega)}
         + \int_{t-1}^{t} (t-s)^{-\frac{1}{2}} \|\boldsymbol{F}(s)\|_{L^q(\varOmega)}\,ds\\
    &\qquad + \int_0^{t-1} (t-s)^{-\frac{3}{2}} \|\boldsymbol{F}(s)\|_{L^q(\varOmega)}\,ds.
  \end{aligned}
  \label{eq:P-22}
 \end{equation}
 Here we have used the standard estimate of analytic semigroup:
 \begin{equation*}
  \|T(t) \vf\|_{W^{1,q}(\varOmega)} \leq C t^{-\frac{1}{2}} \|\vf\|_{L^q(\varOmega)}
 \end{equation*}
 for any $0 < t \leq 1$ and $\vf \in L^q_{\sigma}(\varOmega)$,
 which follows from \eqref{eq:Re-E}.
 By using \eqref{eq:P-17}, \eqref{eq:P-19} and \eqref{eq:P-22} we obtain
 \begin{equation}
  \|\vz(t)\|_{W^{1,q}(\varOmega_{R+1})} \leq C t^{-\frac{3}{2q}}
   \|\vf\|_{L^q(\varOmega)} \quad t \geq 1.
  \label{eq:P-23}
 \end{equation}
 Applying \eqref{eq:P-11} with $r=\infty$ and Lemma~\ref{lem:Bog} we have
 \begin{align*}
  &\|\psi \vw(t)\|_{W^{1,q}(\varOmega_{R+1})} \leq C_q \|\vw(t)\|_{W^{1,\infty}(\R^3)} \leq C_q t^{-\frac{3}{2q}} \| \vf \|_{L^{q}(\varOmega)},\\
  &\begin{aligned}
    \|\B[(\nabla \psi) \cdot \vw(t)]\|_{W^{1,q}(\varOmega_{R+1})}
     &\leq C_q \|(\nabla \psi ) \cdot \vw(t)\|_{L^q(\varOmega_{R+1})}
     \leq C_q \|\vw(t)\|_{L^{\infty}(\R^3)}\\
     &\leq C_q t^{-\frac{3}{2q}} \|\vf\|_{L^q(\varOmega)},
   \end{aligned}
 \end{align*}
 which combined with \eqref{eq:P-12} and \eqref{eq:P-23} implies that
 \begin{equation}
  \|\vv(t)\|_{W^{1,q}(\varOmega_{R+1})} \leq C_q t^{-\frac{3}{2q}} \|\vf\|_{L^q(\varOmega)} \quad \text{for any } t \geq 1.
  \label{eq:P-24}
 \end{equation}
 Now, we shall estimate $\partial_{t} \vv(t)$.
 Recalling that $\vv(t) = T(t+1)\vf \in C^2([0,\infty);\mathcal{D}(\mathcal{M}))$,
 differentiating \eqref{eq:P-6} with respect to $t$ variable,
 we have
 \begin{equation}
 \left\{
	\begin{aligned}
	 &\partial_t \vv_t - \Delta \vv_t = 0, \quad \dv{\vv_t}=0 &\quad &\text{in} &\ &\varOmega \times (0,\infty),\\
	 &\vnu \cdot \vv_t =0, \quad \curl{\vv_t} \times \vnu=0 & &\text{on} &  &\partial\varOmega \times (0,\infty),\\
	 &\vv_t|_{t=0} = {\vg}' & &\text{in} & &\varOmega,
	\end{aligned}
	\right.
	\label{eq:P-25}
 \end{equation}
 where $\vg{}' = \partial_t T(t+1)\vf|_{t=0}$.
 Since $\vg' \in \mathcal{D}(\mathcal{M})$ and $\|\vg'\|_{W^{2,q}(\varOmega)} \leq C_q \| \vf \|_{L^q(\varOmega)}$,
 applying the same argument as above to \eqref{eq:P-25}, we get
 \begin{equation}
  \|\partial_{t} \vv(t)\|_{W^{1,q}(\varOmega_{R+1})}
   \leq C_q t^{-\frac{3}{2q}} \|\vf\|_{L^q(\varOmega)}.
   \label{eq:P-26}
 \end{equation}
 
 Finally we shall estimate the second derivative of $\vv(t)$.
 In order to do this, we shall use Theorem~\ref{thm:resolvent-estimate}
 with $\lambda = 1$.
 Let $\varphi \in C_0^{\infty}(\R^3)$ such that
 $\varphi(x)=1$ for $|x| \leq R$ and $\varphi(x)=0$ for $|x| \geq R+ 1/2$.
 Put 
 \begin{equation*}
  \vv_1(t) = \varphi \vv(t) - \B[(\nabla \varphi) \cdot \vv(t)].
 \end{equation*}
 Here and in the followings, we use the abbreviations $\B \equiv
 \B_{D_{R+1/2,R}}$.
 By Lemma~\ref{lem:Bog} and  Lemma~\ref{prop:d-ext} we have 
 \begin{equation}
  \vv_1(t) = \vv(t) \quad \text{in}\ \varOmega_R,\quad
  \dv{\vv_1(t)} =0 \quad \text{in}\ \varOmega.
  \label{eq:P-27}
 \end{equation}
 According to \eqref{eq:P-6}, \eqref{eq:P-27} and the fact that
 $\vv_1(t) =0$ for $x \not\in B_{R+1/2}$, we have
 \begin{equation*}
 \left\{
	\begin{aligned}
	 &\vv_1(t) - \Delta \vv_1(t) = \boldsymbol{G}(t), \quad \dv{\vv_1}=0 &\quad &\text{in} &\ &\varOmega_{R+1} \times (0,\infty),\\
	 &\vnu \cdot \vv_1(t) =0, \quad \curl{\vv_1(t)} \times \vnu=0 & &\text{on} &  &\partial\varOmega_{R+1} \times (0,\infty),
	\end{aligned}
	\right.
 \end{equation*}
 where 
 \begin{equation*}
  \boldsymbol{G}(t) = \varphi \vv(t) - \B[(\nabla \varphi) \cdot \vv(t)]
   - 2 \nabla \vv(t) \cdot \nabla \varphi - (\Delta \varphi) \vv(t)
   + \Delta \B[(\nabla \varphi) \cdot \vv(t)] + \varphi \partial_{t} \vv(t).
 \end{equation*}
 By Proposition~\ref{prop:BD} we have
 \begin{equation}
  \|\vv_1(t)\|_{W^{2,q}(\varOmega_{R+1})} \leq
  \|\boldsymbol{G}(t)\|_{L^q(\varOmega_{R+1})}.
  \label{eq:P-28}
 \end{equation}
 Applying \eqref{eq:P-24} and \eqref{eq:P-26},
 we have
 \begin{equation*}
  \|\boldsymbol{G}(t)\|_{L^q(\varOmega_{R+1})} \leq C_q t^{-\frac{3}{2q}} \|\vf\|_{L^q(\varOmega)},
 \end{equation*}
 which combined with \eqref{eq:P-27} and \eqref{eq:P-28} implies that
 \begin{equation}
  \|\vv(t)\|_{W^{2,q}(\varOmega_R)} \leq C_{q,R} t^{-\frac{3}{2q}} \|\vf\|_{L^q(\varOmega)}.
  \label{eq:P-29}
 \end{equation}
 Combining \eqref{eq:P-3}, \eqref{eq:P-26} and \eqref{eq:P-29},
 we complete the proof of the lemma.
\end{proof}

\subsection*{2nd step}
At this step, we shall show the following lemma.
\begin{lemma}
 \label{lem:P-2}
 Let $1 < q < \infty$ and $\vf \in L^q_{\sigma}(\varOmega)$.
 Then we have the following two estimates{\normalfont :}
 \begin{equation}
  \|T(t) \vf\|_{L^r(\varOmega)} \leq C_{q,r} t^{-\frac{3}{2}\left(\frac{1}{q} - \frac{1}{r}\right)} \|\vf\|_{L^q(\varOmega)} \quad \text{for any } t \geq 2
  \label{eq:P-30}
 \end{equation}
 provided that $q \leq r \leq \infty$ and $3(1/q - 1/r) <2$ and 
 \begin{equation}
  \|\nabla T(t) \vf\|_{L^q(\varOmega)} \leq C_{q} t^{-\frac{1}{2}} \| \vf\|_{L^q(\varOmega)}\quad \text{for any } t \geq 2
  \label{eq:P-31}
 \end{equation}
 provided that $1 < q \leq 3$.
\end{lemma}

\begin{proof}
 In view of Lemma~\ref{lem:P-1},
 it suffices to estimate $T(t)\vf$ in $\varOmega \setminus B_R$ for $t \geq 2$.
 Set $\vv(t) = T(t+1)\vf=T(t)\vg$ with $\vg = T(1)\vf$.
 Let $\varphi(x) \in C^{\infty}(\R^3)$ so that 
 $\varphi(x)=1$ for $|x| \geq R-1$  and
 $\varphi(x)=0$ for $|x| \leq R-2$.
 In view of Lemma~\ref{prop:d-ext}, we set $\vw(t) = \varphi \vv(t) - \B[(\nabla \varphi) \cdot \vv(t)]$ and then by \eqref{eq:P-6} and Lemma~\ref{lem:Bog} we have
 \begin{equation}
 \left\{
 \begin{aligned}
  &\vw_t - \Delta \vw = \boldsymbol{K}(t), \quad \dv{\vw} = 0 &\quad
  &\text{in} &\ &\R^3 \times (0,\infty),\\ 
  &\vw(0) = \vw_0
 \end{aligned}
 \right.\label{eq:P-32}
 \end{equation}
 where
 \begin{equation}
  \begin{gathered}
   \boldsymbol{K}(t) = - 2 \nabla \vv(t) \cdot \nabla \varphi(t) - (\Delta \varphi) \vv(t) - (\partial_{t} - \Delta) \B[(\nabla \varphi) \cdot \vv(t)],\\
   \vw_0 = \varphi \vg - \B[(\nabla \varphi) \cdot \vg].
  \end{gathered}
  \label{eq:P-33}
 \end{equation}
 Here and hereafter $\B \equiv \B_{R-2,R-1}$. 
 Since $\vw(t) = \vv(t)$ for $|x| \geq R$,
 it suffices to estimate \eqref{eq:P-33}.
 Employing the same arguments as in the proof of \eqref{eq:P-16} and
 \eqref{eq:P-18}, we get
 \begin{gather}
  \dv{\boldsymbol{K}(t)}=0, \dv{\vw_0} =0 \quad \text{in} \ \R^3,
  \label{eq:P-34}\\
  \supp{\boldsymbol{K}(t) \subset D_{R-2,R-1}}.
  \label{eq:P-35}
 \end{gather}
 Let $E(t)$ be the Gaussian kernel: \eqref{eq:Gauss}.
 In view of \eqref{eq:P-34}, 
 employing the same argument as is the proof of \eqref{eq:P-21},
 we have
 \begin{equation}
  \vw(t) = E(t) \ast \vw_0 + \int_0^t E(t-s) \ast \boldsymbol{K}(s)\,ds.
	\label{eq:P-36} 
 \end{equation}
 Applying Young's inequality, we have
 \begin{equation}
  \|\nabla^j E(t) \ast \varphi\|_{L^r(\R^3)} \leq 
  C_{q,r} t^{-\frac{3}{2}\left(\frac{1}{q}-\frac{1}{r}\right)-\frac{j}{2}} \|\varphi\|_{L^q(\R^3)}
  \label{eq:P-37}
 \end{equation}
 for any $t>0$, $j \geq 0$ and $1 \leq q \leq r \leq \infty$,
 and 
 \begin{equation}
 \|E(t) \ast \varphi\|_{W^{2,q}(\R^3)} \leq C t^{-\frac{1}{2}} \|\varphi\|_{W^{1,q}(\R^3)}\label{eq:P-38}
 \end{equation}
 for $0 < t \leq 2$.
 Recalling that $\vv(t) = T(t+1)\vf$, by \eqref{eq:Re-E} we have
 \begin{equation}
  \|\partial_{t} \vv(t)\|_{L^q(\varOmega)} + \|\vv(t)\|_{W^{2,q}(\varOmega)}
  \leq C_q \|\vf\|_{L^q(\varOmega)}
  \label{eq:P-39}
 \end{equation}
 for $0 < t \leq 2$.
 From \eqref{eq:P-4}, \eqref{eq:P-35}, \eqref{eq:P-39},
 Lemma~\ref{lem:Bog} and Lemma~\ref{lem:P-1}, we have
 \begin{gather}
  \|\boldsymbol{K}(t)\|_{W^{1,q}(\R^3)} + \|\boldsymbol{K}(t)\|_{L^{\gamma}(\R^3)} \leq C (1+t)^{-\frac{3}{2q}} \| \vf \|_{L^q(\varOmega)}, \quad 1 \leq \gamma \leq q, \label{eq:P-40}\\
  \|\vw_0\|_{L^q(\R^3)} \leq C \|\vf\|_{L^q(\varOmega)}.
  \label{eq:P-41}
 \end{gather}
 Set
 \begin{equation*}
  I_1(t) = E(t) \ast \vw_0, \quad I_2(t) = \int_0^t E(t-s) \ast \boldsymbol{K}(s)\,ds.
 \end{equation*}
 By \eqref{eq:P-37} and \eqref{eq:P-41} we have
 \begin{equation}
 \begin{aligned}
  \|I_1(t)\|_{L^r(\R^3)} &\leq C_{q,r}  t^{-\frac{3}{2}\left(\frac{1}{q}-\frac{1}{r}\right)} \|\vf\|_{L^q(\varOmega)},\\
  \|\nabla I_1(t)\|_{L^q(\R^3)} &\leq C_q t^{-\frac{1}{2}} \|\vf\|_{L^q(\varOmega)}.
 \end{aligned}\label{eq:P-42}
 \end{equation}
 Let $t \geq 1$ and $r,\gamma$ be numbers such that
 \begin{equation}
  q \leq r \leq \infty, \quad 3 \left(\frac{1}{q} - \frac{1}{r}\right) < 2,
  \quad 1 < \gamma < \min{\left(q,\frac{3}{2}\right)}.
  \label{eq:P-43}
 \end{equation}
 Then by the Sobolev embedding theorem, \eqref{eq:P-37}, \eqref{eq:P-38}
 and \eqref{eq:P-40} we have
 \begin{align}
  \|I_2(t)\|_{L^r(\R^3)}
  &\leq C_{q,r} \int_{t-1}^{t} \|E(t-s) \ast \boldsymbol{K}(s)\|_{W^{2,q}(\R^3)}\,ds \notag \\
  & \qquad + \int_{0}^{t-1} \|E(t-s) \ast \boldsymbol{K}(s)\|_{L^{r}(\R^3)}\,ds
     \notag \\
  &\leq C_{q,r} \int_{t-1}^{t} (t-s)^{-\frac{1}{2}} \|\boldsymbol{K}(s)\|_{W^{1,q}(\R^3)}\,ds \notag\\
  & \qquad  + C_{q,r} \int_0^{t-1} (t-s)^{-\frac{3}{2}\left(\frac{1}{\gamma}-\frac{1}{r}\right)} \|\boldsymbol{K}\|_{L^{\gamma}(\R^3)}\,ds \notag \\
  & \leq C_{q,r} \left\{ t^{-\frac{3}{2q}} + \int_0^{t-1} (t-s)^{-\frac{3}{2}\left(\frac{1}{\gamma}-\frac{1}{r}\right)} (1+s)^{-\frac{3}{2q}} \,ds\right\} \|\vf\|_{L^q(\varOmega)}.\label{eq:P-44}
 \end{align}
 Observe that
 \begin{align*}
  \int_0^{t-1} (t-s)^{-\frac{3}{2}\left(\frac{1}{\gamma}-\frac{1}{r}\right)} (1+s)^{-\frac{3}{2q}}\,ds
  &\leq C_{r,\gamma} \int_0^t (1+t-s)^{-\frac{3}{2}\left(\frac{1}{\gamma}-\frac{1}{r}\right)}(1+s)^{-\frac{3}{2q}}\,ds\\
  &= C_{r,\gamma} \int_0^{t/2} (1+t-s)^{-\frac{3}{2}\left(\frac{1}{\gamma}-\frac{1}{r}\right)} (1+s)^{-\frac{3}{2q}}\,ds\\
  &\qquad    + C_{r,\gamma} \int_0^{t/2} (1+\tau)^{-\frac{3}{2}\left(\frac{1}{\gamma}-\frac{1}{r}\right)} (1+t-\tau)^{-\frac{3}{2q}}\,d\tau,
 \end{align*}
 where we have used the change of variable, $t-s = \tau$ in the second term
 in the last relation.
 When $0 < s < t/2$, $1+t-s \geq 1+s$, we have
 \begin{align*}
  \int_0^{t-1} (t-s)^{-\frac{3}{2}\left(\frac{1}{\gamma}-\frac{1}{r}\right)} (1+s)^{-\frac{3}{2q}}\,ds
  &\leq 2 C_{r,\gamma} \left( 1 + \frac{t}{2}  \right)^{-\frac{3}{2}\left(\frac{1}{q}-\frac{1}{r}\right)} \int_0^{t/2} (1+s)^{-\frac{3}{2\gamma}}\,ds\\
  &\leq 2 C_{q,r} (1+t)^{-\frac{3}{2}\left(\frac{1}{q}-\frac{1}{r}\right)}
 \end{align*}
 because $3\gamma/2 > 1$ holds by \eqref{eq:P-43},
 which combined with \eqref{eq:P-44} implies that 
 \begin{equation}
  \|I_2(t)\|_{L^r(\R^3)} \leq C_{q,r} t^{-\frac{3}{2}\left(\frac{1}{q}-\frac{1}{r}\right)} \|\vf\|_{L^q(\varOmega)} \quad 
  \text{for any } t \geq 1
  \label{eq:P-45}
 \end{equation}
 provided that $q \leq r \leq \infty$ and $3(1/q - 1/r) < 2$.
 From \eqref{eq:P-41} and \eqref{eq:P-42} we have
 \begin{equation}
  \|\nabla I_1(t)\|_{L^q(\R^3)} 
  \leq C_q t^{-\frac{1}{2}} \|\vf\|_{L^q(\varOmega)}.
  \label{eq:P-46}
 \end{equation}
 By \eqref{eq:P-37} and \eqref{eq:P-40} we have
 \begin{equation}
  \begin{aligned}
   \|\nabla I_2(t)\|_{L^q(\R^3)} 
   &\leq C_q \int_{t-1}^{t} (t-s)^{-\frac{1}{2}} (1+s)^{-\frac{3}{2q}}\,ds \|\vf\|_{L^q(\varOmega)}\\
   &\qquad + C_{q,r} \int_0^{t-1}
   (t-s)^{-\frac{3}{2}\left(\frac{1}{\gamma}-\frac{1}{q}\right) - \frac{1}{2}} (1+s)^{-\frac{3}{2q}}\,ds \|\vf\|_{L^q(\varOmega)}.
   \end{aligned}
   \label{eq:P-47}
 \end{equation}
 Observe that
 \begin{align*}
   \int_0^{t-1}
  (t-s)^{-\frac{3}{2}\left(\frac{1}{\gamma}-\frac{1}{q}\right) -
  \frac{1}{2}} (1+s)^{-\frac{3}{2q}}\,ds 
   &\leq C_{q,\gamma} \int_0^t (1+t-s)^{-\frac{3}{2}\left(\frac{1}{\gamma}-\frac{1}{q}\right)-\frac{1}{2}} (1+s)^{-\frac{3}{2q}}\,ds \\
   &= C_{q,\gamma} \int_0^{t/2} (1+t-s)^{-\frac{3}{2}\left(\frac{1}{\gamma}-\frac{1}{q}\right)-\frac{1}{2}} (1+s)^{-\frac{3}{2q}}\,ds\\
   & \qquad + C_{q,\gamma} \int_{t/2}^{t} (1+s)^{-\frac{3}{2}\left(\frac{1}{\gamma}-\frac{1}{q}\right)-\frac{1}{2}} (1+t-s)^{-\frac{3}{2q}}\,ds.
 \end{align*}
 If $1 < q \leq 3$, then $3/2q - 1/2 \geq 0$, and therefore
 \begin{align*}
  \int_{0}^{t-1} (t-s)^{-\frac{3}{2}\left(\frac{1}{\gamma}-\frac{1}{q}\right)-\frac{1}{2}} (1+s)^{-\frac{3}{2q}}\,ds
  &\leq C_{q,\gamma}\left(1 + \frac{t}{2}\right)^{-\frac{1}{2}} \int_0^{t/2} (1+s)^{-\frac{3}{2\gamma}}\,ds\\
  &\leq C_{q,\gamma} (1+t)^{-\frac{1}{2}},
 \end{align*}
 which combined with \eqref{eq:P-47} implies that
 \begin{equation}
  \|\nabla I_2(t)\|_{L^q(\R^3)} \leq C_q t^{-\frac{1}{2}} \|\vf\|_{L^q(\varOmega)}, \quad t \geq 1
 \end{equation}
 provided that $1 < q \leq 3$.
 The proof is completed.
\end{proof}

\subsection*{3rd step}
We consider the case when $0 < t \leq 2$.
We shall prove the following lemma.
\begin{lemma}
 \label{lem:P-3}
 Let $1 < q < \infty$ and $0 < t \leq 2$ and $\vf \in L^q_{\sigma}(\varOmega)$.
 Then we have
 \begin{align}
  \|T(t)\vf\|_{L^r(\varOmega)} &\leq C_{q,r} t^{-\frac{3}{2}\left(\frac{1}{q}-\frac{1}{r}\right)} \|\vf\|_{L^q(\varOmega)}, &\quad &1 < q \leq r \leq \infty,\label{eq:P-48}\\
  \|\nabla T(t)\vf\|_{L^r(\varOmega)} &\leq C_{q,r} t^{-\frac{3}{2}\left(\frac{1}{q}-\frac{1}{r}\right)-\frac{1}{2}} \|\vf\|_{L^q(\varOmega)}, &\quad &1 < q \leq r < \infty.\label{eq:P-49}
 \end{align}
\end{lemma}
\begin{proof}
 For any real number $s \in (0,2)$, by  complex interpolation theorem
 we have $W^{s,q}(\varOmega) = [L^q(\varOmega),W^{2,q}(\varOmega)]_{\theta}$ with $s=2\theta$ (see e.g., Triebel \cite{Tr}).
 From \eqref{eq:Re-E} we have
 \begin{align}
  \|T(t)\vf\|_{L^q(\varOmega)} &\leq C_q \|\vf\|_{L^q(\varOmega)},
     \label{eq:P-50}\\
  \|T(t)\vf\|_{W^{2,q}(\varOmega)} &\leq C_q t^{-1} \| \vf \|_{L^q(\varOmega)}
     \label{eq:P-50'}
 \end{align}
 for $0 < t \leq 2$.
 Therefore interpolating \eqref{eq:P-50} and \eqref{eq:P-50'} for $s=2\theta$
 we obtain
 \begin{equation}
  \|T(t)\vf\|_{W^{s,q}(\varOmega)} 
  \leq C_{q,s} t^{-\frac{s}{2}} \|\vf\|_{L^q(\varOmega)}.
  \label{eq:P-52}
 \end{equation}
 From the Sobolev embedding theorem and \eqref{eq:P-52},
 for $s=3(1/q-1/r)$ we have
 \begin{equation}
  \|T(t) \vf \|_{L^r(\varOmega)} \leq C_{q,r} t^{-\frac{3}{2}\left(\frac{1}{q} - \frac{1}{r}\right)} \| \vf \|_{L^q(\varOmega)}
  \label{eq:P-53}
 \end{equation}
 for $0 < t \leq 2$ and $1 < q \leq r < \infty$.
 By \eqref{eq:P-50} and \eqref{eq:P-50'} we have
 \begin{equation}
  \|\nabla T(t) \vf\|_{L^q(\varOmega)} \leq C \|T(t) \vf\|_{L^q(\varOmega)}^{\frac{1}{2}} \|T(t) \vf\|_{W^{2,q}(\varOmega)}^{\frac{1}{2}} 
  \leq C t^{-\frac{1}{2}} \|\vf\|_{L^q(\varOmega)}
  \label{eq:P-54}
 \end{equation}
 for $0 < t \leq 2$.
 Therefore, by \eqref{eq:P-53} and \eqref{eq:P-54} we obtain
 \begin{align*}
  \|\nabla T(t) \vf \|_{L^r(\varOmega)} &\leq \left\|\nabla T \left( \frac{t}{2} + \frac{t}{2} \right) \vf \right\|_{L^r(\varOmega)}
  \leq C \left(\frac{t}{2}\right)^{-\frac{1}{2}} \left\| T\left(\frac{t}{2}\right) \vf \right\|_{L^r(\varOmega)}\\
  &\leq C_{q,r} t^{-\frac{3}{2}\left(\frac{1}{q} - \frac{1}{r}\right) - \frac{1}{2}} \|\vf \|_{L^q(\varOmega)}
 \end{align*}
 for $0 < t \leq 2$.

 Finally we shall consider the $L^{\infty}$ estimate. 
 For $3 < q < \infty$, by using Sobolev's inequality:
 \begin{equation*}
  \|\vu\|_{L^{\infty}(\varOmega)}
  \leq C \|\vu\|_{W^{1,q}(\varOmega)}^{\theta}
  \|\vu\|_{L^q(\varOmega)}^{1-\theta} 
 \end{equation*} 
 with $\theta = 3/q$ and \eqref{eq:P-53} and \eqref{eq:P-54} we have
 \begin{equation}
  \|T(t)\vf\|_{L^{\infty}(\varOmega)} \leq C_q t^{-\frac{3}{2q}} \|\vf\|_{L^q(\varOmega)}
  \label{eq:P-55}
 \end{equation}
 for $0 < t \leq 2$.
 Next we consider the cases when $1 < q < 3/2$ or $3/2 < q < 3$. 
 Let $3/(k+1) < q < 3/k$ with $k=1,2$.
 We set $\{q_\ell\}_{\ell=0}^{k}$ in such a way that $1/{q_{\ell+1}} =
 1/q_{\ell} - 1/3\ (\ell=0,1,\dots,k-1)$ with $q_0 = q$.
 Since $1 < q < 3$, we see that $3 < q_k < \infty$.
 Therefore by using \eqref{eq:P-55} with $q = q_{k}$ and \eqref{eq:P-53}
 with $r = q_{k}$, we obtain
 \begin{align*}
  \|T(t)\vf\|_{L^{\infty}(\Omega)} &=
  \left\|T\left(\frac{t}{2}\right)T\left(\frac{t}{2}\right)
  \vf\right\|_{L^{\infty}(\Omega)} \leq C t^{-\frac{3}{2 q_{k}}}
  \left\|T\left(\frac{t}{2}\right)\right\|_{L^{q_k}(\Omega)}\\
  &\leq C t^{-\frac{3}{2 q_k}} t^{-\frac{3}{2} \left(\frac{1}{q} -
  \frac{1}{q_k}\right)} \|\vf\|_{L^q(\Omega)} = C t^{-\frac{3}{2q}}
  \|\vf\|_{L^q(\Omega)},
 \end{align*}
 for $t > 0$.
 This implies \eqref{eq:P-55} for $1 - 3/q \not\in \mathbb{N}_0$.
 When $1 - 3/q \in \mathbb{N}_0$, 
 we choose $r$ in such a way that $q < r < \infty$ and $1 - 3/r \not\in \mathbb{N}_0$.
 Then, by \eqref{eq:P-53} with $q=r$ and \eqref{eq:P-55} we have
 \begin{align*}
  \|T(t) \vf\|_{L^{\infty}(\varOmega)}
  &\leq  C_r \left( \frac{t}{2}\right)^{-\frac{3}{2r}} \left\| T\left(\frac{t}{2}\right) \vf\right\|_{L^r(\varOmega)} \\
  &\leq C_r \left(\frac{t}{2}\right)^{-\frac{3}{2r}} \left(\frac{t}{2}\right)^{-\frac{3}{2} \left(\frac{1}{q} - \frac{1}{r}\right)} \|\vf\|_{L^q(\varOmega)}  \leq C_{q,r} t^{-\frac{3}{2q}} \|\vf\|_{L^q(\varOmega)}
 \end{align*}
 for $0 < t \leq 2$.
 Hence we get \eqref{eq:P-53} for $1 < q \leq r \leq \infty$.
 The proof is completed.
\end{proof}

\subsection*{4th step}
Now, we shall complete the proof of Theorem~\ref{thm:lq-lr}.
Combining Lemma~\ref{lem:P-2} and Lemma~\ref{lem:P-3}, 
we have 
\begin{equation}
 \|T(t) \vf \|_{L^r(\varOmega)} \leq
 C_{q,r} t^{-\frac{3}{2}\left(\frac{1}{q} - \frac{1}{r}\right)} \|\vf\|_{L^q(\varOmega)}
 \label{eq:P-56}
\end{equation}
for any $t > 0$ and $\vf \in L^q_{\sigma}(\varOmega)$ 
provided that $1 < q \leq r \leq \infty$ and $3(1/q-1/r)<2$.
When $1 < q \leq r \leq \infty$ and $3(1/q-1/r) \geq 2$,
we choose numbers $q_j$, $j=1,2,\dots,\ell-1$,
in such a way that $q=q_0 < q_1 < q_2 < \dots< q_{\ell-1} < q_{\ell} = r$
and $3(1/q_{m-1} - 1/q_{m}) < 2$ for $m=1,2,\dots,\ell$.
Repeated use of \eqref{eq:P-56} implies that
\begin{align*}
  \|T(t) \vf\|_{L^r(\varOmega)} &=  \biggl\|T\biggl(\underbrace{\frac{t}{\ell}+ \cdots + \frac{t}{\ell}}_{ \text{$\ell$ times}}\biggr)\vf \biggr\|_{L^r(\varOmega)} \\
  &\leq C_{q_{\ell},q_{\ell-1}} \left(\frac{t}{\ell}\right)^{-\frac{3}{2}\left( \frac{1}{q_{\ell-1}} - \frac{1}{r}\right)} \biggl\|T\biggl(\underbrace{\frac{t}{\ell}+ \cdots + \frac{t}{\ell}}_{ \text{$\ell-1$ times}}\biggr)\vf \biggr\|_{L^{q_{\ell-1}}(\varOmega)} \\
  &\leq \dots \leq C_{q,r} \left(\frac{t}{\ell}\right)^{-\frac{3}{2}\left(\frac{1}{q} - \frac{1}{r}\right)} \|\vf\|_{L^q(\varOmega)},
\end{align*}
and therefore we have \eqref{eq:P-56} 
for any $t > 0$ and $\vf \in L^q_{\sigma}(\varOmega)$
provided that $1 < q \leq r \leq \infty$.

Now we consider the case when $q=1$.
For any 
$\boldsymbol{\varphi},\boldsymbol{\psi} \in C_{0,\sigma}^{\infty}(\varOmega)$,
by Theorem~\ref{thm:B-B*} and \eqref{eq:P-56}
we have
\begin{equation*}
 |(T(t)\boldsymbol{\varphi},\boldsymbol{\psi})_{\varOmega}| = |(\boldsymbol{\varphi},T(t)\boldsymbol{\psi})_{\varOmega}| \leq \|\boldsymbol{\varphi}\|_{L^1(\varOmega)} \|T(t)
  \boldsymbol{\psi}\|_{L^{\infty}(\varOmega)} \leq C \|\boldsymbol{\varphi}\|_{L^1(\varOmega)} t^{-3/2r'} \|\boldsymbol{\psi}\|_{L^{r'}(\varOmega)},
\end{equation*}
where $r'=r/(r-1)$,
and therefore we have
\begin{equation}
 \|T(t) \boldsymbol{\varphi}\|_{L^r(\varOmega)} \leq C_r t^{-\frac{3}{2}\left(1-\frac{1}{r}\right)} \|\boldsymbol{\varphi}\|_{L^1(\varOmega)}.
 \label{eq:P-57}
\end{equation}
Since $C_{0,\sigma}^{\infty}(\varOmega)$ is dense in $L^{r'}_{\sigma}(\varOmega)$,
by the density argument we have \eqref{eq:P-57} for any
$\boldsymbol{\varphi} \in L^1_{\sigma}(\varOmega) = \overline{C_{0,\sigma}^{\infty}(\Omega)}^{\|\cdot\|_{L^1(\Omega)}}$.

Combining \eqref{eq:P-31} and \eqref{eq:P-49}, we obtain
\begin{equation}
 \|\nabla T(t) \vf\|_{L^q(\varOmega)} \leq C_q t^{-\frac{1}{2}} \| \vf \|_{L^q(\varOmega)}
 \label{eq:P-58}
\end{equation}
for any $t > 0$ and $\vf \in L^q_{\sigma}(\varOmega)$
provided that $1 < q \leq 3$.
Combining \eqref{eq:P-56}, \eqref{eq:P-57} and \eqref{eq:P-58},
we have
\begin{align*}
 \|\nabla T(t) \vf\|_{L^r(\varOmega)}
 &= \left\| \nabla T
 \left(\frac{t}{2}\right) T \left(\frac{t}{2}\right) \vf 
  \right\|_{L^r(\varOmega)}
 \leq C_r t^{-\frac{1}{2}} \left\|T \left(\frac{t}{2}\right)
   \vf \right\|_{L^r(\varOmega)}\\
 &\leq C_{q,r} t^{-\frac{3}{2}\left(\frac{1}{q} - \frac{1}{r}\right) - \frac{1}{2}} \| \vf \|_{L^{q}(\varOmega)}.
\end{align*}
for any $t > 0$ and $\vf \in L^q_{\sigma}(\varOmega)$
provided that $1 \leq q \leq r \leq 3$, $r \not= 1$.
This completes the proof of Theorem~\ref{thm:lq-lr}.

\section{Proof of Theorem~\ref{thm:global}}
\label{sec:global-existence}
This section is devoted to the proof of Theorem~\ref{thm:global}.
For notational simplicity,
we use the abbreviation $\|\cdot\|_{q}$ which stands for
$\|\cdot\|_{L^q(\varOmega)}$.
At first employing the argument due to Kato \cite{Kato-84} for the Cauchy
problem of the Navier-Stokes system,
 we shall solve the integral equations \eqref{eq:INT} 
by contraction mapping principle.

In order to do this, we introduce the following symbols:
\begin{align*}
 [\vv]_{\ell,q,t} &= \sup_{0 < s \leq t} s^{\ell} \|\vv(s)\|_{q},\\
 \LN \vv \RN_{t} &= [\vv]_{\frac{1-\delta}{2},\frac{3}{\delta},t} + [\nabla \vv]_{\frac{1}{2},3,t},\\
 \+ \vv \+_{t} &= [\vv]_{0,3,t} + [\vv]_{\frac{1}{2},\infty,t}+\LN \vv
 \RN_{t} 
\end{align*}
with some fixed real number $\delta \in (0,1)$.
As an underlying space, we set
\begin{align}
 \mathcal{I}_{M} = \{ &(\vv(t),\vB(t)) \in BC([0,\infty); L^3_{\sigma}(\varOmega) \times L^3_{\sigma}(\varOmega)) \,|\,\notag \\
 & \quad \lim_{t \rightarrow 0+}\{ [(\vv - \va,\vB-\vb)]_{0,3,t} +
 [(\vv,\vB)]_{\frac{1}{2},\infty,t}+\LN  (\vv,\vB) \RN_{t}
 \}=0,\label{eq:I-1} \\  
 & \quad \sup_{t>0}\+(\vv,\vB)\+_{t} \leq 2 M
 \|(\va,\vb)\|_{3}\},\label{eq:I-2} 
\end{align}
where $M$ will be determined later (see \eqref{eq:4.8} below).
Set 
\begin{align*}
 &\vv_0(t) =e^{-tA}\va, \quad \vB_0(t) = e^{-t\mathcal{M}}\vb,\\
 &\Phi(\vv,\vB)(t)=\binom{\vv_0(t)}{\vB_0(t)} + \binom{F[\vv,\vB](t)}{G[\vv,\vB](t)}.
\end{align*}
We shall prove that there exist positive constants $M$ and $\eta$ such
that if
\begin{equation}
 \|(\va,\vb)\|_{3} \leq \eta,\label{eq:smallness}
\end{equation}
then $\Phi$ becomes a contraction map from $\mathcal{I}_M$ 
into itself.

At the beginning, we shall show that
\begin{gather}
 \lim_{t \rightarrow 0+} [(\vv_0-\va,\vB_0-\vb)]_{0,3,t}=0,
  \label{eq:C00-1}\\
 \lim_{t \rightarrow 0+} \LN (\vv_0,\vB_0) \RN_{t}=0, \qquad
  \lim_{t \rightarrow 0+} [ (\vv_0,\vB_0) ]_{\frac{1}{2},\infty,t}=0.
  \label{eq:C00-2}
\end{gather}
In fact, for any $\epsilon > 0$ there exists a pair
$(\va_{\epsilon},\vb_{\epsilon}) \in C_{0,\sigma}^{\infty}(\varOmega)
\times C_{0,\sigma}^{\infty}(\Omega)$ so that 
$\|(\va,\vb)-(\va_\epsilon,\vb_\epsilon)\|_{3} < \epsilon$.
Therefore, by the $L^3$-boundedness of the semigroups
(Theorems~\ref{thm:lq-lr} and \ref{thm:lq-lr-Stokes} with $q=r=3$), we
see that
\begin{align*}
 \|(\vv_0(t),\vB_0(t)) - (\va,\vb)\|_{3}
 &\leq \|(e^{-tA}(\va-\va_\epsilon),e^{-t\mathcal{M}}(\vb-\vb_\epsilon))\|_{3}\\
 &\qquad+ \|(e^{-tA}\va_\epsilon-\va_\epsilon, e^{-t\mathcal{M}}\vb_\epsilon - \vb)\|_{3}
        + \|(\va_\epsilon-\va,\vb_\epsilon-\vb)\|_{3}\\
 &\leq C\epsilon 
        + \|(e^{-tA}\va_\epsilon-\va_\epsilon,
 e^{-t\mathcal{M}}\vb_\epsilon - \vb)\|_{3}\\ 
 &\leq C\epsilon + \int_0^t \left\| \frac{d}{ds} (e^{-sA}\va_\epsilon,e^{-s\mathcal{M}}\vb_\epsilon)  \right\|_{3}\,ds\\
 &\leq C \epsilon + C t \|(\va_\epsilon,\vb_\epsilon)\|_{W^{2,3}(\varOmega)}.
\end{align*}
Therefore we have
\begin{equation*}
 \lim_{t \rightarrow 0+}  [(\vv_0-\va,\vB_0-\vb)]_{0,3,t} \leq C \epsilon.
\end{equation*}
This implies \eqref{eq:C00-1}, because $\epsilon$ is chosen arbitrarily.
By similar manner, we have
\begin{align*}
 t^{\frac{1-\delta}{2}} \|(\vv_0(t),\vB_0(t))\|_{\frac{3}{\delta}}
 &\leq t^{\frac{1-\delta}{2}} \|(e^{-tA}(\va-\va_\epsilon),e^{-t\mathcal{M}}(\vb-\vb_\epsilon))\|_{\frac{3}{\delta}} \\
 &\qquad + t^{\frac{1-\delta}{2}} \|(e^{-tA}\va_\epsilon,e^{-t\mathcal{M}} \vb_\epsilon)\|_{\frac{3}{\delta}}\\
 &\leq C \|(\va-\va_\epsilon,\vb-\vb_\epsilon)\|_{3}
      + C t^{\frac{1}{2}-\frac{3}{2r}} \|(\va_\epsilon,\vb_\epsilon)\|_{r}\\
 &\leq C \epsilon + C t^{\frac{1}{2}-\frac{3}{2r}} \|(\va_\epsilon,\vb_\epsilon)\|_{r}
\end{align*}
with some $r \in (3,3/\delta)$,
which implies 
\begin{equation}
 \lim_{t \rightarrow 0+} [\vv_0,\vB_0]_{\frac{1-\delta}{2},\frac{3}{\delta},t} \leq C \epsilon. \label{eq:4.6}
\end{equation}
From similar calculation, we see that 
 \begin{equation}
  \lim_{t \rightarrow 0+} [(\vv_0,\vB_0)]_{\frac{1}{2},\infty,t}  \leq 
  C \epsilon,\qquad 
  \lim_{t \rightarrow 0+} [\nabla (\vv_0,\vB_0)]_{\frac{1}{2},3,t} \leq
  C \epsilon. \label{eq:4.7}
 \end{equation}
Since $\epsilon$ is chosen arbitrarily,
by \eqref{eq:4.6} and \eqref{eq:4.7} we have
\eqref{eq:C00-2}. 

By Theorems~\ref{thm:lq-lr} and \ref{thm:lq-lr-Stokes},
one can easily see that
\begin{align}
 \+ (\vv_0,\vB_0) \+_{t} \leq M \|(\va,\vb)\|_{3} \quad \text{for any } t>0
 \label{eq:4.8}
\end{align}
with some constant $M$.
In particular,  from \eqref{eq:4.6}, \eqref{eq:4.7} and
\eqref{eq:4.8},
 we see that $(\vv_0(t),\vB_0(t)) \in \mathcal{I}_M$.

Now, we shall estimate the nonlinear terms
$F[\vv,\vB](t)$ and $G[\vv,\vB](t)$.
In order to do this, we prepare the following inequality essentially due
to the H\"older inequality:
\begin{equation}
 \|(\vu(s) \cdot \nabla)\vv(s)\|_{\frac{3}{1+\delta}}
 \leq \|\vu(s)\|_{\frac{3}{\delta}} \|\nabla \vv(s)\|_{3}
 \leq C s^{-1 + \frac{\delta}{2}} \LN \vu \RN_t \LN \vv \RN_t
\label{eq:4.9}
\end{equation}
for any $0 < s \leq t$.
By Theorem~\ref{thm:lq-lr-Stokes} and the $L^q$-boundedness of the
Helmholtz projection \eqref{eq:HP-Boundedness},
we have
\begin{align*}
 \|F[\vv,\vB](t)\|_{3} 
 &\leq \int_0^t \|e^{-(t-s)A} P[(\vv(s) \cdot \nabla)\vv(s) - (\vB(s) \cdot \nabla)\vB(s)]\|_{3}\,ds\\
 &\leq C \int_{0}^{t} (t-s)^{-\frac{3}{2}\left(\frac{1+\delta}{3} - \frac{1}{3}\right)} (\|(\vv(s) \cdot \nabla)\vv(s)\|_{\frac{3}{1+\delta}} 
        + \|(\vB(s) \cdot \nabla)\vB(s)\|_{\frac{3}{1+\delta}})\,ds\\
 &\leq C \int_0^t (t-s)^{-\frac{\delta}{2}} (\|\vv(s)\|_{\frac{3}{\delta}} \|\nabla \vv(s)\|_{3} + \|\vB(s)\|_{\frac{3}{\delta}}\|\nabla \vB(s)\|_{3})\,ds.
\end{align*}
By similar manner with Theorem~\ref{thm:lq-lr}, we have
\begin{align*}
 \|G[\vv,\vB](t)\|_{3} 
 &\leq C \int_0^t (t-s)^{-\frac{\delta}{2}} (\|\vv(s)\|_{\frac{3}{\delta}} \|\nabla \vB(s)\|_{3} + \|\vB(s)\|_{\frac{3}{\delta}}\|\nabla \vv(s)\|_{3})\,ds.
\end{align*}
From the above two estimates and \eqref{eq:4.9}, we obtain
\begin{align}
 \|(F[\vv,\vB])(t),G[\vv,\vB](t)\|_{3}
 &\leq C \int_0^t (t-s)^{-\frac{\delta}{2}} \|(\vv(s),\vB(s))\|_{\frac{3}{\delta}} \|\nabla (\vv(s),\vB(s))\|_{3}\,ds\notag\\
 &\leq C \int_0^t (t-s)^{-\frac{\delta}{2}} s^{-1+\frac{\delta}{2}}\,ds
 \LN (\vv,\vB)\RN^2_{t}\notag\\ 
 &= C B\left(1-\frac{\delta}{2},\frac{\delta}{2}\right) \LN
 (\vv,\vB)\RN_t^2,
 \label{eq:4.a}
\end{align}
where $B(q,r)$ denotes the beta function.
From similar calculations, we obtain the following estimates:
\begin{align}
 \|(F[\vv,\vB](t),G[\vv,\vB](t)\|_{\frac{3}{\delta}}
 &\leq C \int_0^t (t-s)^{-\frac{1}{2}} s^{-1+\frac{\delta}{2}}\,ds \LN
 (\vv,\vB)\RN_t^2\notag \\ 
 &\leq C B\left(\frac{1}{2},\frac{\delta}{2}\right) t^{-\frac{1-\delta}{2}} \LN (\vv,\vB)\RN_t^2;\label{eq:4.b}\\
 \|\nabla (F[\vv,\vB](t),G[\vv,\vB](t)\|_{3}
 &\leq C \int_0^t (t-s)^{-\frac{1+\delta}{2}} s^{-1+\frac{\delta}{2}}\,ds \LN (\vv,\vB)\RN_t^2\notag\\
 &\leq C B\left(\frac{1-\delta}{2},\frac{\delta}{2}\right)
 t^{-\frac{1}{2}} \LN (\vv,\vB)\RN_t^2;\label{eq:4.c}\\
 \|(F[\vv,\vB](t),G[\vv,\vB](t)\|_{\infty} 
 &\leq C \int_0^t (t-s)^{-\frac{1+\delta}{2}} s^{-1+\frac{\delta}{2}}\,ds \LN (\vv,\vB)\RN_t^2\notag\\
 &\leq C B\left(\frac{1-\delta}{2},\frac{\delta}{2}\right) t^{-\frac{1}{2}} \LN (\vv,\vB)\RN_t^2.\label{eq:4.d}
\end{align}
From \eqref{eq:4.a}, \eqref{eq:4.b}, \eqref{eq:4.c} and \eqref{eq:4.d},
we have 
\begin{equation}
 \+ (F[\vv,\vB],G[\vv,\vB]) \+_{t} \leq C \LN (\vv,\vB)\RN_t^2.
 \label{eq:4.10}
\end{equation}
Hence, from \eqref{eq:4.8} and \eqref{eq:4.10},
we have
\begin{align}
 &\+ \Phi(\vv,\vB)\+_{t} \leq M \|(\va,\vb)\|_{3} + C \LN (\vv,\vB)\RN_t^2,
 \label{eq:4.12}\\
 &\begin{aligned} 
  &[\Phi(\vv,\vB)-(\va,\vb)] {}_{0,3,t} +
   [\Phi(\vv,\vB)]_{\frac{1}{2},\infty,t} + \LN \Phi(\vv,\vB)\RN_t\\
  \leq &\ 
  [(\vv_0-\va,\vB_0-\vb)]_{0,3,t} + [(\vv_0,\vB_0)]_{\frac{1}{2},\infty,t}
   + \LN (\vv_0,\vB_0)\RN_{t} + C \LN (\vv,\vB)\RN_t^2.
   \end{aligned}
\label{eq:4.13}
\end{align}
Therefore, if $(\vv,\vB) \in \mathcal{I}_M$,
then by \eqref{eq:I-1}, \eqref{eq:I-2}, \eqref{eq:C00-1}, \eqref{eq:C00-2},
\eqref{eq:4.12} and \eqref{eq:4.13}, we obtain
\begin{gather}
 \+ \Phi(\vv,\vB) \+_{t} \leq M \|(\va,\vb)\|_{3} + 4CM^2
 \|(\va,\vb)\|_3^2 \quad \text{for any } t>0,
 \label{eq:4.14}\\
 \lim_{t \rightarrow 0+} ([\Phi(\vv,\vB)-(\va,\vb)]_{0,3,t} +
 [\Phi(\vv,\vB)]_{\frac{1}{2},\infty,t}+  \LN
 \Phi(\vv,\vB)\RN_t)=0.\label{eq:4.15} 
\end{gather}
Choose an $\eta > 0$ in such a way that
\begin{equation}
 4 C M \eta < 1 \label{eq:4.16}.
\end{equation}
Then by \eqref{eq:4.14} we have
\begin{equation}
 \+ \Phi(\vv,\vB) \+_{t} < 2 M \|(\va,\vb)\|_{3} \quad \text{for any } t>0
 \label{eq:4.16'}
\end{equation}
provided that $\|(\va,\vb)\|_{3} \leq \eta$,
which combined with \eqref{eq:4.15} implies that
$\Phi(\vv,\vB) \in \mathcal{I}_M$ provided that $(\vv,\vB) \in \mathcal{I}_M$.
This shows that $\Phi$ is a mapping from $\mathcal{I}_M$ into itself.
By using \eqref{eq:4.9}, Theorems~\ref{thm:lq-lr} and \ref{thm:lq-lr-Stokes}
and employing the same argument as in the proof
of \eqref{eq:4.10}, we have
\begin{equation}
\begin{aligned}
 &\+ \Phi(\vv_1,\vB_1) - \Phi(\vv_2,\vB_2) \+_t \\
 \leq \,&C (\LN(\vv_1,\vB_1)\RN_{t}  + \LN(\vv_2,\vB_2)\RN_{t}) \LN(\vv_1,\vB_1)-(\vv_2,\vB_2)\RN_{t} \\
 \leq \,&4 CM \|(\va,\vb)\|_{3} \+(\vv_1,\vB_1)-(\vv_2,\vB_2)\+_{t}
\end{aligned}
\label{eq:4.17}
\end{equation}
for any $(\vv_1,\vB_1), (\vv_2,\vB_2) \in \mathcal{I}_M$.
If we choose an $\eta>0$ in such a way that
\begin{equation*}
 4 C M \eta < \frac{1}{2},
\end{equation*}
then it follows from \eqref{eq:4.17} that $\Phi$ is a 
contraction map from $\mathcal{I}_M$ into itself
if $\|(\va,\vb)\|_{3} \leq \eta$.
Therefore, there exists 
a unique fixed point $(\vv(t),\vB(t)) \in \mathcal{I}_M$
of $\Phi$, which solves \eqref{eq:INT}.
The uniqueness of solutions to \eqref{eq:INT} holds 
for any $(\vv(t),\vB(t)) \in \mathcal{I}_M$.
Namely, if $(\vv_1(t),\vB_1(t)), (\vv_2(t),\vB_2(t)) \in \mathcal{I}_M$
satisfy the integral equations \eqref{eq:INT} with
the same initial data $(\va,\vb) \in L^3_{\sigma}(\varOmega) \times L^3_{\sigma}(\varOmega)$ with $\|(\va,\vb)\|_{3} \leq \eta$,
then we have $(\vv_1(t),\vB_1(t)) = (\vv_2(t),\vB_2(t))$ for any $t>0$.

Now we shall show sharp asymptotic behavior of the global in time strong
solution: \eqref{eq:AS-1} and \eqref{eq:AS-2}.
In order to do this, at first we shall show the following:
\begin{equation}
 \lim_{t \rightarrow \infty} \|(\vv(t),\vB(t))\|_{3} =0.
\label{eq:G-5}
\end{equation}
Given $ 0 < \gamma < 1/2$, we take $3/2 < q < 3$ such that
$\gamma = 3/2q - 1/2$.
Given $(\va,\vb) \in C^{\infty}_{0,\sigma}(\varOmega) \times C_{0,\sigma}^{\infty}(\Omega)$ with
$\|(\va,\vb)\|_{3} < \eta$, let $(\vv(t),\vB(t))$ be
 solution of \eqref{eq:INT}.
Then applying the $L^q$-$L^3$ estimate and the $L^{3/2}$-$L^3$ estimate for
$e^{-tA}$ and $e^{-t\mathcal{M}}$ to \eqref{eq:INT} and using the H\"older
inequality, we have
\begin{equation*} 
\begin{aligned}
 &\|(\vv(t),\vB(t))\|_{3} \\
 &\ \leq C t^{-\gamma} \|(\va,\vb)\|_{q} + C \int_0^t (t-s)^{-\frac{1}{2}} \|(\vv(s),\vB(s))\|_{3} \|\nabla (\vv(s),\vB(s))\|_{3}\,ds\\
 &\ \leq C t^{-\gamma} \|(\va,\vb)\|_{q} 
  +   C \int_0^t (t-s)^{-\frac{1}{2}} s^{-\gamma} s^{-\frac{1}{2}}\,ds
   [(\vv,\vB)]_{\gamma,3,t} 
   [\nabla (\vv,\vB)]_{\frac{1}{2},3,t}\\
 &\ \leq C t^{-\gamma} \left\{ \|(\va,\vb)\|_{q} 
    + C B\left(\frac{1}{2},\frac{1}{2}-\gamma\right) \|(\va,\vb)\|_{3}
 \|[(\vv,\vB)]_{\gamma,3,t} \right\},
\end{aligned}
\end{equation*}
which implies that 
\begin{equation*}
 [(\vv,\vB)]_{\gamma,3,t} \leq C \|(\va,\vb)\|_{q} + C \|(\va,\vb)\|_{3}[(\vv,\vB)]_{\gamma,3,t}.
\end{equation*}
Since choosing $\eta>0$ smaller if necessary,
we may assume that $C\|(\va,\vb)\|_{3} \leq 1/2$
provided that $\|(\va,\vb)\|_{3} \leq \eta$,
we have
\begin{equation*}
 [(\vv,\vB)]_{\gamma,3,t} \leq 2 C \|(\va,\vb)\|_{q}.
\end{equation*}
This implies that \eqref{eq:G-5} holds for any 
initial data  $(\va,\vb) \in C_{0,\sigma}^{\infty}(\varOmega) \times
C_{0,\sigma}^{\infty}(\Omega)$ with $\|(\va,\vb)\|_{3} \leq \eta$.

For general $(\va,\vb) \in L^3_{\sigma}(\varOmega) \times
L^3_{\sigma}(\Omega)$ with 
$\|(\va,\vb)\|_{3} < \eta$ and any $\epsilon > 0$,
we choose $\va_{\epsilon}$ and $\vb_{\epsilon}$ in such a way that
$\|(\va_{\epsilon}-\va,\vb_{\epsilon}-\vb)\|_{3} \leq \epsilon$.
Choosing $\epsilon>0$ smaller if necessary,
we may assume that
$\|(\va_{\epsilon},\vb_{\epsilon})\|_{3} < \eta$ for any $\epsilon > 0$.
Since $\|(\va_{\epsilon},\vb_{\epsilon})\|_{3} < \eta$,
the corresponding solution of \eqref{eq:INT} satisfies \eqref{eq:G-5}.
Combining this fact and continuous dependence of solution:
$L^3_{\sigma}(\varOmega) \times L^3_{\sigma}(\Omega) \ni (\va,\vb)
\mapsto (\vv(t),\vB(t)) \in BC([0,\infty);L^3_{\sigma}(\varOmega) \times
L^3_{\sigma}(\Omega))$,
we have 
\begin{align*}
 \|(\vv(t),\vB(t))\|_{3} &\leq \|(\vv(t)-\vv_{\epsilon}(t),\vB(t)-\vB_{\epsilon}(t))\|_3 + \|(\vv_{\epsilon}(t),\vB_{\epsilon}(t))\|_3\\
 &\leq C \epsilon + C \|(\vv_\epsilon(t),\vB_\epsilon(t))\|_{3}.
\end{align*}
Since $\epsilon$ is arbitrary and 
 $(\vv_\epsilon(t),\vB_\epsilon(t))$ satisfies \eqref{eq:G-5},
we get \eqref{eq:G-5} for any initial data 
$(\va,\vb) \in L^3_{\sigma}(\varOmega) \times L^3_{\sigma}(\Omega)$ with
 $\|(\va,\vb)\|_{3} \leq \eta$.

By the interpolation inequality, we get
\begin{align*}
 t^{\frac{1}{2}-\frac{3}{2q}} \|(\vv(t),\vB(t))\|_{q}
 &\leq \|(\vv(t),\vB(t))\|_{3}^{\theta} \left( t^{\frac{1}{2}} \|(\vv(t),\vB(t))\|_{\infty} \right)^{1-\theta}\\
  &\leq C_q \|(\va,\vb)\|^{1-\theta}_{3} \|(\vv(t),\vB(t))\|_3^{\theta}
\end{align*}
with $1/q=\theta/3$,
which together with \eqref{eq:G-5} implies that
\eqref{eq:AS-1} for $3 < q < \infty$.
Here we have used the global boundedness of $t^{1/2}
\|(\vv(t),\vB(t))\|_{\infty}$ which is guaranteed by the fact that a
pair $(\vv(t),\vB(t))$ is global solution of \eqref{eq:INT} with
property \eqref{eq:I-1} and \eqref{eq:I-2}. 
Finally, we shall prove 
\eqref{eq:AS-1} for $q = \infty$ and \eqref{eq:AS-2}.
In order to do this, we rewrite \eqref{eq:INT} as follows:
\begin{equation*}
\left\{
\begin{aligned}
  \vv(t) &= e^{-\frac{t}{2}A} \vv (t/2)- \int_{\frac{t}{2}}^{t} e^{-(t-s)A} P[(\vv(s) \cdot \nabla)\vv(s) - (\vB(s) \cdot \nabla) \vB(s)]\,ds,\\
 \vB(t) &= e^{-\frac{t}{2}\mathcal{M}} \vB (t/2)- \int_{\frac{t}{2}}^{t} 
 e^{-(t-s)\mathcal{M}} [(\vv(s) \cdot \nabla)\vB(s) - (\vB(s) \cdot
 \nabla) \vv(s)]\,ds. 
\end{aligned}
\right.
\end{equation*}
Then by Theorems~\ref{thm:lq-lr} and \ref{thm:lq-lr-Stokes} we obtain
\begin{align*}
\|(\vv(t),\vB(t))\|_{\infty}
\leq \,&C t^{-\frac{1}{2}} \left\| \left( \vv(t/2), \vB (t/2)  \right)\right\|_{3}\\
 &+ C \int_{\frac{t}{2}}^{t} (t-s)^{-\frac{3}{4}} \|(\vv(s),\vB(s))\|_{6} \|\nabla (\vv(s),\vB(s))\|_{3}\,ds
\end{align*}
and 
\begin{align*}
\|\nabla(\vv(t),\vB(t))\|_{3}
\leq \,&C t^{-\frac{1}{2}} \left\|\left( \vv(t/2), \vB (t/2)  \right)\right\|_{3}\\
 &+ C \int_{\frac{t}{2}}^{t} (t-s)^{-\frac{3}{4}} \|(\vv(s),\vB(s))\|_{6} \|\nabla (\vv(s),\vB(s))\|_{3}\,ds
\end{align*}
Therefore combining the above two estimates and $\|\nabla (\vv(t),\vB(t))\|_{3} \leq C t^{-1/2} \|(\va,\vb)\|_{3}$,
we obtain
\begin{equation*}
\begin{aligned}
& t^{\frac{1}{2}} (\|(\vv(t),\vB(t))\|_{\infty} + \|\nabla (\vv(t),\vB(t))\|_{3})\\
 &\ \leq C \|(\vv(t/2),\vB(t/2))\|_{3} + C \|(\va,\vb)\|_{3} \sup_{t/2 \leq s \leq t} s^{\frac{1}{4}} \|(\vv(s),\vB(s))\|_{6}
\end{aligned}
\end{equation*}
for $t > 0$.
Therefore, from \eqref{eq:G-5} and \eqref{eq:AS-1} with $q=6$,
we have \eqref{eq:AS-1} for $q = \infty$ and \eqref{eq:AS-2}.
This completes the proof of Theorem~\ref{thm:global}.

\bigskip
\noindent
\textit{Acknowledgments}.\quad
The author would like to express his heartily gratitude to Professor
Yoshihiro Shibata for valuable comments and constant encouragements.
The author is also grateful to the referee for many valuable comments
and helpful suggestions.


\begin{thebibliography}{10}

\bibitem{AKST-04}
T.~Akiyama, H.~Kasai, Y.~Shibata, and M.~Tsutsumi.
\newblock On a resolvent estimate of a system of {L}aplace operators with
  perfect wall condition.
\newblock {\em Funkcial. Ekvac.}, 47(3):361--394, 2004.

\bibitem{Bg}
M.~E. Bogovski{\u\i}.
\newblock Solution of the first boundary value problem for an equation of
  continuity of an incompressible medium.
\newblock {\em Dokl. Akad. Nauk SSSR}, 248(5):1037--1040, 1979.

\bibitem{B-S}
W.~Borchers and H.~Sohr.
\newblock On the semigroup of the {S}tokes operator for exterior domains in
  ${L}\sp q$-spaces.
\newblock {\em Math. Z.}, 196(3):415--425, 1987.

\bibitem{D-L-1972}
G.~Duvaut and J.-L. Lions.
\newblock In\'equations en thermo\'elasticit\'e et magn\'eto-hydrodynamique.
\newblock {\em Arch. Rational Mech. Anal.}, 46:241--279, 1972.

\bibitem{Enomoto-Shibata-05-JMFM}
Y.~Enomoto and Y.~Shibata.
\newblock On the rate of decay of the {O}seen semigroup in exterior domains and
  its application to {N}avier-{S}tokes equations.
\newblock {\em J. Math. Fluid Mech.}, 7(3):339--367, 2005.

\bibitem{GlI}
G.~P. Galdi.
\newblock {\em An introduction to the mathematical theory of the
  {N}avier-{S}tokes equations. {V}ol. {I}}.
\newblock Springer-Verlag, New York, 1994.
\newblock Linearized steady problems.

\bibitem{G-M-85}
Y.~Giga and T.~Miyakawa.
\newblock Solutions in {$L\sb r$} of the {N}avier-{S}tokes initial value
  problem.
\newblock {\em Arch. Rational Mech. Anal.}, 89(3):267--281, 1985.

\bibitem{MR991022}
Y.~Giga and H.~Sohr.
\newblock On the {S}tokes operator in exterior domains.
\newblock {\em J. Fac. Sci. Univ. Tokyo Sect. IA Math.}, 36(1):103--130, 1989.

\bibitem{MR2085848}
T.~Hishida.
\newblock The nonstationary {S}tokes and {N}avier-{S}tokes flows through an
  aperture.
\newblock In {\em Contributions to current challenges in mathematical fluid
  mechanics}, Adv. Math. Fluid Mech., pages 79--123. Birkh\"auser, Basel, 2004.

\bibitem{Iwashita}
H.~Iwashita.
\newblock ${L}\sb q$-${L}\sb r$ estimates for solutions of the nonstationary
  {S}tokes equations in an exterior domain and the {N}avier-{S}tokes initial
  value problems in ${L}\sb q$ spaces.
\newblock {\em Math. Ann.}, 285(2):265--288, 1989.

\bibitem{Kato-84}
T.~Kato.
\newblock Strong {${L}\sp{p}$}-solutions of the {N}avier-{S}tokes equation in
  {${\R}^{m}$}, with applications to weak solutions.
\newblock {\em Math. Z.}, 187(4):471--480, 1984.

\bibitem{Kozono-87}
H.~Kozono.
\newblock On the energy decay of a weak solution of the {M}{H}{D} equations in
  a three-dimensional exterior domain.
\newblock {\em Hokkaido Math. J.}, 16(2):151--166, 1987.

\bibitem{L-S-1960}
O.~A. Ladyzhenskaya and V.~A. Solonnikov.
\newblock Solution of some non-stationary problems of magnetohydrodynamics for
  a viscous incompressible fluid.
\newblock {\em Trudy Mat. Inst. Steklov}, 59:115--173, 1960.

\bibitem{L-L}
L.~D. Landau and E.~M. Lifshitz.
\newblock {\em Electrodynamics of continuous media}.
\newblock Course of Theoretical Physics, Vol. 8. Translated from the Russian by
  J. B. Sykes and J. S. Bell. Pergamon Press, Oxford, 1960.

\bibitem{M-S}
P.~Maremonti and V.~A. Solonnikov.
\newblock On nonstationary {S}tokes problem in exterior domains.
\newblock {\em Ann. Scuola Norm. Sup. Pisa Cl. Sci. (4)}, 24(3):395--449, 1997.

\bibitem{My}
T.~Miyakawa.
\newblock The {$L\sp{p}$} approach to the {N}avier-{S}tokes equations with the
  {N}eumann boundary condition.
\newblock {\em Hiroshima Math. J.}, 10(3):517--537, 1980.

\bibitem{My2}
T.~Miyakawa.
\newblock On nonstationary solutions of the {N}avier-{S}tokes equations in an
  exterior domain.
\newblock {\em Hiroshima Math. J.}, 12(1):115--140, 1982.

\bibitem{Pazy}
A.~Pazy.
\newblock {\em Semigroups of linear operators and applications to partial
  differential equations}, volume~44 of {\em Applied Mathematical Sciences}.
\newblock Springer-Verlag, New York, 1983.

\bibitem{S-T}
M.~Sermange and R.~Temam.
\newblock Some mathematical questions related to the {M}{H}{D} equations.
\newblock {\em Comm. Pure Appl. Math.}, 36(5):635--664, 1983.

\bibitem{MHD-local-energy}
Y.~Shibata and N.~Yamaguchi.
\newblock Local energy decay for a parabolic system related to Maxwell's
		equations in exterior domain.
\newblock \textit{preprint}.

\bibitem{MR1190728}
C.~G. Simader and H.~Sohr.
\newblock A new approach to the {H}elmholtz decomposition and the {N}eumann
  problem in {$L\sp q$}-spaces for bounded and exterior domains.
\newblock In {\em Mathematical problems relating to the Navier-Stokes
  equation}, volume~11 of {\em Ser. Adv. Math. Appl. Sci.}, pages 1--35. World
  Sci. Publishing, River Edge, NJ, 1992.

\bibitem{Tr}
H.~Triebel.
\newblock {\em Interpolation theory, function spaces, differential operators}.
\newblock Johann Ambrosius Barth, Heidelberg, second edition, 1995.

\bibitem{5th-RIMS}
N.~Yamaguchi.
\newblock ${L}^q$-${L}^r$ estimates of solution to the parabolic {M}axwell
  equations and their application to the magnetohydrodynamic equations.
\newblock {\em S\=urikaisekikenky\=usho K\=oky\=uroku}, (1353):72--91, 2004.
\newblock Mathematical analysis in fluid and gas dynamics (Japanese) (Kyoto,
  2003).

\bibitem{Y-G}
Z.~Yoshida and Y.~Giga.
\newblock On the {O}hm-{N}avier-{S}tokes system in magnetohydrodynamics.
\newblock {\em J. Math. Phys.}, 24(12):2860--2864, 1983.

\end{thebibliography}
\end{document}